\documentclass{article}

\usepackage{amsmath}
\usepackage{amscd}
\usepackage{amssymb}

\usepackage{amsthm}
    
\usepackage[usenames]{color}    
\newtheorem{Def}{Definition}[section]     
\newtheorem{Prop}[Def]{Proposition}
\newtheorem{Lemma}[Def]{Lemma}     
\newtheorem{Thm}[Def]{Theorem} 
\newtheorem{Cor}[Def]{Corollary} 
\newtheorem{Rem}[Def]{Remark}

\newcommand{\C}{\mathbb{C}}
\newcommand{\R}{\mathbb{R}}
\newcommand{\Z}{\mathbb{Z}}

\newcommand{\rad}{\mathrm{rad}}
\newcommand{\Ee}{\tilde{\mathbb E}}
\newcommand{\Hh}{\mathbb{H}}
\newcommand{\Pp}{\mathbb{P}}
\newcommand{\PD}{\Pp(D^2)}
\newcommand{\PH}{\Pp(\Hh_{half})}

\newcommand{\Dkai}{\stackrel{\circ}{D^2}}

\newcommand{\Specan}{\mathrm{Specan}}

\newcommand{\de}{{\delta}}

\newcommand{\Cc}{\mathcal{C}}
\newcommand{\Aa}{\mathcal{A}}

\newcommand{\Oo}{\mathcal {O}}

\newcommand{\Ff}{\mathcal {F}}
\newcommand{\Ss}{\hbox{\boldmath$\mathcal {S}$}}

\newcommand{\bp}{\begin{pmatrix}}
\newcommand{\ep}{\end{pmatrix}}

\numberwithin{equation}{section}

\newcounter{CounterEQUlabel}
\newcommand{\EQUlabel}[1]{\label{#1}
	\ifcase \theCounterEQUlabel
		\relax
	\or
		\hspace{1em}\mbox{\small$\langle$\rmfamily#1$\rangle$}
		\index{zzz#1@#1}
	\fi }	
	\newcounter{CounterEQUref}
	\newcounter{CounterEQUpageref}
	\newcommand{\EQUref}[1]{
		\ifcase \theCounterEQUref     \relax   \or {\small[#1]}\,\fi
		\ifcase \theCounterEQUpageref (\ref{#1})\or (\ref{#1}\,(p.\pageref{#1})) \fi}

\newcounter{Counterakakaigyou}
\newcommand{\akakaigyou}[1]{
	\ifcase \theCounterakakaigyou
		\relax
	\or
		{\par
		\color{red}
		 #1 
		\color{black}
		\par}
	\fi }

\newcounter{Counteraka} 
\newcommand{\aka}[1]{
	\ifcase \theCounteraka
		\relax
	\or
		{
		\color{red}
		 #1 
		\color{black}
		}
	\fi }

\newcounter{Counteraokaigyou}
\newcommand{\aokaigyou}[1]{
	\ifcase \theCounteraokaigyou
		\relax
	\or
		{\par
		\color{blue}
		 #1 
		\color{black}
		\par}
	\fi }	

\newcounter{Counterao} 
\newcommand{\ao}[1]{
	\ifcase \theCounterao
		\relax
	\or
		{
		\color{blue}
	     #1 
		\color{black}
		}
	\fi }	
\newcounter{CounterEQUnewpage} 
\newcommand{\EQUnewpage}{
	\ifcase \theCounterEQUnewpage
		\relax
	\or
		{\newpage}
	\fi }		

\newcommand{\EQUspace}{
	\ifcase \theCounterEQUnewpage
		\relax
	\or
		{\vskip1cm}
	\fi }	
	
\title{An $O(2,n)$ formulation of invariant theory for elliptic Weyl group}
\author{Ikuo Satake}
\date{Department of Mathematics, Osaka University, 
Toyonaka, Osaka, 560-0043, JAPAN}
\begin{document}
\maketitle
\setcounter{section}{0}
\setcounter{secnumdepth}{3}
\setcounter{tocdepth}{5}

\setcounter{CounterEQUlabel}{0} 
\setcounter{CounterEQUref}{0}
\setcounter{CounterEQUpageref}{0}

\setcounter{Counterakakaigyou}{0} 
\setcounter{Counteraka}{0} 
\setcounter{Counteraokaigyou}{0} 
\setcounter{Counterao}{0} 
\setcounter{CounterEQUnewpage}{0}

\begin{abstract}
In this paper, we give a new formulation of 
invariant theory for elliptic Weyl group 
using the group $O(2,n)$. 
As an elliptic Weyl group quotient, 
we define a suitable $\C^*$-bundle. 
We show that it has a conformal Frobenius 
structure which we define in this paper. 
Then its good section could be identified 
with a Frobenius manifold which 
we constructed in \cite{SatakeFrobenius}. 
\end{abstract}

\section{Introduction}
Motivated by a period mapping for a primitive form 
\cite{Saito}
for an unfolding of 
a function with a simple elliptic singularity, 
an elliptic root system is introduced \cite{extendedI}. 
On the quotient space $\Ee_a//W$ of the domain $\Ee_a$ by an 
elliptic Weyl group $W$, a flat holomorphic metric is 
constructed by \cite{extendedII} and 
Frobenius manifold structure is constructed by \cite{SatakeFrobenius}. 

Here the domain $\Ee_a$ is not canonical, 
that is, it depends on the choice of a marking $a$. 
Then we set the following problems.

Problem 1:
What is a relation between the elliptic Weyl group quotient spaces 
$\Ee_a//W$ and $\Ee_{a'}//W^2$ for different markings $a$ and $a'$?

Problem 2:
What is a relation between the Frobenius manifold structure on them?

By definition, an elliptic root system is a root system belonging 
to a vector space with an inner produce with a signature 
$(0,2,l)$ where we denote by $(l_+,l_0,l_-)$ 
numbers of positive, $0$ and negative eigenvalues respectively
(we choose a negative semi-definite inner product). 
A domain $\Ee_a$ corresponds to a hyperbolic extension with a 
signature $(1,1,l+1)$. 

In this paper, we consider 2-extension with a signature 
$(2,0,l+2)$ and construct a $\C^*$-bundle 
\begin{equation}
p:D^2 \to \PD,
\end{equation}
where $D^2$ is a quadric hypersurface of a suitable domain $D$ 
and $\PD=D^2/\C^*$. They corresponds to the group $O(2,0,l+2)$. 
Then we find that $\Ee_a$ could be canonically identified 
with a section of $p$. 

This formulation is good from the viewpoint of 
the group action of the automorphism of an elliptic root system. 
We find that a central extension group of 
the group of the automorphisms of an elliptic root system 
is naturally defined in our formulation and this group 
acts on the $\C^*$-bundle $p$ equivariantly. 
The fact that this action changes a section of $p$, 
corresponds to the fact that this group is not 
contained in $O(1,1,l+1)$. 

We see a Frobenius manifold structure. 
The elliptic Weyl group is constructed as a subgroup 
of $O(2,0,n)$ and denote it by $W^2$, which is isomorphic to $W$. 
Then the Weyl group quotient:
\begin{equation}
p:D^2//W^2 \to \PD//W^2
\end{equation}
is again a $\C^*$-bundle and the space $\Ee_a//W$ gives 
a section of it. 
Since $\Ee_a//W$ is a Frobenius manifold, 
the base space $\PD//W^2$ is also a Frobenius manifold. 
But if we choose another marking $a'$ and 
a corresponding another section $\Ee_{a'}//W$, 
then $\PD//W^2$ has a structure of another Frobenius manifold. 
In fact, except holomorphic metrics, these two Frobenius manifold 
structures are same. For holomorphic metrics, 
they are conformally equivalent. 

Thus we introduce a notion of conformal Frobenius structure 
and we prove that $p:D^2//W^2 \to \PD//W^2$ has 
essentially unique conformal Frobenius structure. 
For any marking $a$, a Frobenius manifold structure of $\Ee_a//W$ 
is obtained from this conformal Frobenius manifold structure. 
This clarifies the relation of the Frobenius manifold 
structures of $\Ee_a//W$ and $\Ee_{a'}//W^2$. 

Also we see that the conformal Frobenius structure has 
an action of a (central extension of) automorphism group 
of the elliptic root system. 
This clarifies the meaning of a conformal transformation 
in \cite{Satakeauto} (see Remark \ref{5.2remark2}).

This paper is organaized as follows. 

In Section 2, we define 2-extension for elliptic root system and 
define an elliptic Weyl group $W^2$ 
and a $\C^*$-bundle $p:D^2 \to \PD$. 
We introduce open subsets of a $\C^*$-bundle 
$p:\Dkai \to \Pp(\Dkai)$   
whose complements are zero of roots. 
Then we take an elliptic Weyl group 
quotient $p:\Dkai/W^2 \to \Pp(\Dkai)/W^2$. 
We call them ``analytic Weyl group quotient spaces''. 
On the space $\Dkai/W^2$, we define a tensor $I^*_{\Dkai/W^2}$ 
derived from an inner product of the elliptic root system. 

In Section 3, we define elliptic Weyl group invariant rings. 
We define ``algebraic Weyl group quotient spaces'' 
$p:D^2//W^2 \to \PD//W^2$ by taking $\Specan$ 
of these elliptic Weyl group invariant rings. 
On the space $D^2//W^2$, we also define a tensor $I^*_{D^2//W^2}$. 
derived from an inner product of the elliptic root system. 
Our theorems are stated for ``algebraic Weyl group quotient spaces''. 
But for proofs, we utilize ``analytic Weyl group quotient spaces''. 

In Section 4, we introduce a notion of conformal Frobenius structure 
and state that a $\C^*$-bundle $p:D^2//W^2 \to \PD//W^2$ 
has essentially unique conformal Frobenius structure. 
We also show that a group of central extension of automorphism 
group of the elliptic root system acts on the conformal 
Frobenius structure. 

In Section 5, We give proofs. We relate our formulation with 
a formulation of 
\cite{extendedII} and \cite{SatakeFrobenius} 
in which the space $\Ee_a$ is used. 
We obtain the existence of a 
conformal Frobenius structure by this comparison. 
For a uniqueness theorem, we use a conformal deformation 
of a holomorphic metric of a Frobenius manifold. 

We use the following terminologies. 

For a complex manifold $M$, 
we denote by $\Oo_M$ 
(resp. $\Omega^1_M,\ \Theta_M$)
the sheaf of holomorphic functions 
(resp. holomorphic 1-forms, holomorphic vector fields). 
We denote 
by $\Oo(M)$ (resp. $\Omega^1(M),\ \Theta(M)$) 
the module $\Gamma(M,\Oo_M)$ 
(resp. $\Gamma(M,\Omega^1_M),\ \Gamma(M,\Theta_M)$).

Let $\Ff$ be a sheaf on the complex manifold $N$. 
For a morphism $f:M \to N$ of complex manifolds 
$M$ and $N$, we denote by $f^*$ the following 
natural morphism 
\begin{equation}\EQUlabel{1.1equation}
f^*:\Gamma(N,\Ff) \to \Gamma(M,f^*\Ff).
\end{equation}
induced by $\Ff \to f_*f^*\Ff$. 

For an isomorphism $f:M \to M$, we have canonical isomorphisms 
$f^*\Omega^1_M \stackrel{\sim}{\longrightarrow}
\Omega^1_M$ and 
$f^*\Theta_M \stackrel{\sim}{\longleftarrow}
\Theta_M$. Then the morphism \EQUref{1.1equation} induces 
\begin{equation}
f^*:\Gamma(M,
(\Omega^1_M)^{\otimes p}
\otimes_{\Oo_M}
(\Theta_M)^{\otimes q}
) 
\stackrel{\sim}{\longrightarrow}
\Gamma(M,
(\Omega^1_M)^{\otimes p}
\otimes_{\Oo_M}
(\Theta_M)^{\otimes q}
).
\end{equation}

For a $\C^*$-bundle $p:L \to M$ and the isomorphism $f$ 
of the $\C^*$-bundle:
\begin{equation}
\begin{CD}
L @>{f}>{\sim}> L\\
@V{p}VV @V{p}VV\\
M @>{f}>{\sim}> M
\end{CD},
\end{equation}
we have 
$f^*(p^*\Omega^1_M) \simeq p^*f^*\Omega^1_M 
\stackrel{\sim}{\longrightarrow}
p^*\Omega^1_M$ 
and 
$f^*(p^*\Theta_M) \simeq p^*f^*\Theta_M 
\stackrel{\sim}{\longleftarrow}
p^*\Theta_M$. 
Then the morphism \EQUref{1.1equation} induces 
\begin{equation}
f^*:\Gamma(L,
(p^*\Omega^1_M)^{\otimes p}
\otimes_{\Oo_L}
(p^*\Theta_M)^{\otimes q}
) 
\stackrel{\sim}{\longrightarrow}
\Gamma(L,
(p^*\Omega^1_M)^{\otimes p}
\otimes_{\Oo_L}
(p^*\Theta_M)^{\otimes q}
).
\end{equation}

\section{Analytic Weyl group quotient space}

\subsection{Elliptic root system}

In this subsection, we define an elliptic root system, its orientation and 
2-extension.

\subsubsection{Definition of elliptic root system}
Let $l$ be a positive integer. 
Let $F$ be a real vector space of rank $l+2$ 
with a negative semi-definite or 
positive semi-definite 
symmetric bilinear form 
$I_F:F \times F \to \R$, whose radical 
$\mathrm{rad}I_F:=\{x \in F\,|\,I_F(x,y)=0,\forall y \in F\}$ 
is a vector space of rank 2. 
For a non-isotropic element $\alpha \in F$ (i.e. 
$I_F(\alpha,\alpha) \neq 0$), we put 
$\alpha^{\vee}:=2\alpha/I_F(\alpha,\alpha) \in F$. 
The reflection $w_{\alpha}$ with respect to $\alpha$ is defined by 
\begin{equation}
w_{\alpha}(u):=u-I_F(u,\alpha^{\vee})\alpha\quad
(\forall u \in F).
\end{equation}
%
%
\begin{Def}
{\rm 
(\cite[p.104, Def. 1]{extendedI}) }
A set $R$ of non-isotropic elements of $F$ is an elliptic 
root system belonging to $(F,I_F)$ if it satisfies the axioms 
1-4:
\begin{enumerate}
	\item The additive group generated by $R$ in $F$, 
	denoted by $Q(R)$, is a full sub-lattice of $F$. 
	That is, the embedding $Q(R) \subset F$ induces 
	the isomorphism : $Q(R) \otimes_{\Z}\R \simeq 
	F$. 
	\item $I_F(\alpha,\beta^{\vee}) \in \Z$ for 
	$\alpha,\beta \in R$.
	\item $w_{\alpha}(R)=R$ for $\forall \alpha \in R$. 
	\item If $R=R_1 \cup R_2$, with $R_1 \perp R_2$, 
	then either $R_1$ or $R_2$ is void. 
\end{enumerate}
\end{Def}
\subsubsection{Definition of orientation}
For an elliptic root system $R$ belonging to $(F,I_F)$, 
the additive group $\mathrm{rad}I_F \cap Q(R)$ 
is isomorphic to $\Z^2$. 
\begin{Def}
An elliptic root system $R$ is called oriented 
if the $\R$-vector space $\mathrm{rad}I_F$ is oriented. 
An ordered pair $\{a,b\}$ of $\mathrm{rad}I_F$ is called 
an oriented basis if it gives an $\R$-basis of $\mathrm{rad}I_F$ 
and it gives the orientation of $\mathrm{rad}I_F$. 
\end{Def}

\subsubsection{Definition of a signed marking}
\begin{Def}
Let $R$ be an elliptic root system $R$ belonging to $(F,I_F)$. 
By a signed marking, we mean a non-zero element $a$ of 
$\mathrm{rad}I_F \cap Q(R)$ such that $Q(R) \cap \R a=\Z a$ 
and the quotient root system 
$R/\R a$ ($:=\mathrm{Image}(R \hookrightarrow 
F \to F/\R a$)) is reduced (i.e. $\alpha,c\alpha \in R/\R a$ 
implies $c \in \{\pm 1\}$). 
We denote the set of a signed marking by $\Lambda_{\Z}$. 
\end{Def}
Hereafter we assume that $\Lambda_{\Z} \neq \emptyset$. 
\EQUspace
\subsubsection{2-extension}
We define a 2-extension. 
Let $F^2$ be a real vector space of rank $l+4$ and 
$I_{F^2}:F^2 \times F^2 \to \R$ an $\R$-symmetric bilinear form. 
The pair $(F^2,I_{F^2})$ is called a 2-extension of 
$(F,I_F)$ if 
$F^2$ contains $F$ as a linear subspace, 
$\mathrm{rad}I_{F^2}=\{0\}$ and $I_{F^2}|_{F}=I_F$. 
$I_{F^2}$ is indefinite symmetric bilinear form 
with signature of $(2,l+2)$ or $(l+2,2)$ according to 
$I$ is negative semi-definite or positive semi-definite. 
A 2-extension is unique up to isomorphism. 
Hereafter we fix a hyperbolic extension $(F^2,I_{F^2})$.
\EQUnewpage
\subsection{Orthogonal group and Weyl group}
\subsubsection{Definition of an Orthogonal group and a Weyl group}
In this subsection, we define a parabolic subgroup of 
an orthogonal group and Weyl group.

Let $F^2_{\C}:=F^2 \otimes_{\R}\C$, 
$(\rad I_F)_{\C}:=\rad I_F\otimes_{\R}\C$. 
$I_{F^2_{\C}}$ is a $\C$-bilinear extension of $I_{F^2}$. 
We put $\tilde{w}_{\alpha}(u):=u-I_{F^2_{\C}}(u,\alpha^{\vee})\alpha$ 
for $u \in F^2_{\C}$. 

We define the groups
\begin{align}
SL(\rad I_F)&:=\{g \in GL(\rad I_F)\,|\,\det g=1\,\},\\
O(F)&:=\{g \in GL(F)\,|\,I_F(gx,gy)=I_F(x,y),\,\forall 
x,y \in F\,\},\\
O^+(F)&:=\{g \in O(F)\,|\,g|_{\rad I_F} \in SL(\rad I_F)\,\},\\
Aut(R)&:=\{g \in GL(F)\,|\,g(R)=R,\ g \in O^+(F)\},\\
%
%
%
%
O(F^2_{\C})&:=\{g \in GL(F^2_{\C})\,|\,
I_{F^2_{\C}}(gx,gy)=I_{F^2_{\C}}(x,y),\ 
\forall x,y \in F^2_{\C}\},\\
O^+(F^2_{\C},F)&:=\{g \in O(F^2_{\C})\,|\,
g(F) \subset F,\ g|_F \in O^+(F)\},\\
W^2&:=
\langle \tilde{w}_{\alpha} \in O^+(F^2_{\C},F)\,|\,
\alpha \in R \rangle.
\end{align}

The natural morphism
\begin{equation}
\pi^2:O^+(F^2_{\C},F) \to O^+(F)\ 
\end{equation}
is surjective. 
We put 
\begin{align}
K^2_{\C}&:=\mathrm{ker}\pi^2,\\
K^2_{\Z}&:=\mathrm{ker}(\pi^2|_{W^2}),\\
\widetilde{Aut}(R)&:=(\pi^2)^{-1}(Aut(R)).
\end{align}

\subsubsection{Properties of an orthogonal group}
\begin{Prop}
1) The group $O^+(F^2_{\C},F)$ acts on the flag 
$F^2_{\C} \supset F \supset \rad I_F$. \\
2) The group $K^2_{\C}$ is isomorphic to additive group $\C$ 
and it is contained in the center of $O^+(F^2_{\C},F)$. \\
3) The group 
$\widetilde{Aut}(R)$ contains $W^2$ as a normal subgroup. \\
\end{Prop}
\begin{proof}
By definition, we obtain 1). 
We obtain 2) by direct calculation. 
Since 
$g\tilde{w}_{\alpha}g^{-1}
=\tilde{w}_{g(\alpha)}$ for 
$g \in \widetilde{Aut}(R)$, we obtain 3).
\end{proof}

\EQUnewpage
\subsection{Description of a center}

In this subsection, we construct the canonical group 
isomorphism $\rho:\C \to K_{\C}$ by the aid of 
the Eichler-Siegel transformation (\cite[p93]{extendedI}). 

We define the Eichler-Siegel transformation by
\begin{equation}
ES:F^2_{\C} \otimes_{\C}F^2_{\C} \to \mathrm{End}_{\C}(F^2_{\C}),\quad
\sum_i \alpha_i \otimes \beta_i 
\mapsto 
u-\sum_i \alpha_iI_{F^2_{\C}}(u,\beta_i).
\end{equation}
Here $F^2_{\C}\otimes_{\C}F^2_{\C}$ has 
the semi-group structure by the product
\begin{equation}
(\sum_{i}u_i \otimes v_i)
\circ
(\sum_{j}w_j \otimes x_j)
=
(\sum_{i}u_i \otimes v_i)
+
(\sum_{j}w_j \otimes x_j)
-
(\sum_{i,j}I_{F^2_{\C}}(v_i,w_j)u_i \otimes x_j).
\end{equation}
Then $ES$ is a semi-group isomorphism.

\begin{Prop}\EQUlabel{2.3prop}
1) $ES$ induces the isomorphism 
$\wedge^2(\rad I_F)_{\C}$ to 
$K^2_{\C}$. 
\\
2) A semigroup
\begin{equation}
\R_{>0}(\mathrm{sgn}(I_F))(a \otimes b-b \otimes a)
\cap ES^{-1}(W^2)\,(\simeq \Z_{>0})
\end{equation}
has the unique generator $g_0$, where $\mathrm{sgn}(I_F)$ is $\pm 1$ 
according to $I_F$ is positive semi-definite or negative semi-definite.\\
3) The morphism $ES$ and $g_0$ depend on the normalization of $I_F$. 
But the morphism 
\begin{equation}
\rho:\C \to K_{\C},\quad
\alpha \mapsto ES(g_0 \otimes \alpha)
\end{equation}
does not depend on a scalar multiplication 
($\R^{\times}$) of $I_{F^2}$. 
\end{Prop}
\begin{proof}
We could interprete the morphism 
$\pi^2:O^+(F^2_{\C},F) \to O^+(F)$ 
into the spaces $F^2_{\C}\otimes_{\C}F^2_{\C}$ 
and $F \otimes_{\R}F/\rad I_F$ 
by the morphism $ES$. 
Then we have a result. 
\end{proof}

\begin{Rem}\EQUlabel{2.3remark}
The element $g_0$ is exactly written by the notion 
of marked extended affine root system as 
$g_0=(I_R:I)\frac{l_{max}+1}{m_{max}}(a \otimes b-b \otimes a)$, 
if $a,b$ are oriented basis 
and $\mathrm{rad}I_F \cap Q(R) \simeq \Z a \oplus \Z b$. 
See \cite[p27 (2.6.2)]{extendedII}.
\end{Rem}

\EQUnewpage
\subsection{$\C^*$-bundles}
\subsubsection{Preparations}
We arbitrary fix an oriented basis $a,b$. 
We put 
$(F^2_{\C})^*:=\mathrm{Hom}_{\C}(F^2_{\C},\C)$, 
$(\rad I_F)^*_{\C}:=\mathrm{Hom}_{\C}((\rad I_F)_{\C},\C)$. 
We denote a natural pairing $F^2_{\C} \times (F^2_{\C})^* \to \C$ 
by $\langle \ ,\ \rangle$. 
A non-degenerate $\C$-blinear form $I_{(F^2_{\C})}$ gives 
a dual $\C$-blinear form 
$I_{(F^2_{\C})^*}:(F^2_{\C})^* \times (F^2_{\C})^* \to \C$. 

The group $GL(F^2_{\C})$ gives the left action on the space 
$(F^2_{\C})^*$ by 
\begin{equation}
\langle
g^{-1}\cdot \alpha,x
\rangle
=
\langle
\alpha,g\cdot x
\rangle\EQUlabel{2.4action}
\end{equation}
for 
$
\forall \alpha \in F^2_{\C},
\forall x \in (F^2_{\C})^*,
\forall g \in GL(F^2_{\C})
$. 

\EQUspace

\subsubsection{Definition of $\C^*$ bundles}

We put
\begin{align}
&D:=\{x \in (F^2_{\C})^*\,|\,
	\langle a,x \rangle \neq 0,\,
	\langle b,x \rangle \neq 0,\,
	\mathrm{Im}\frac{\langle b,x \rangle}
	{\langle a,x \rangle}>0\,\},\\
&D^2:=\{x \in D\,|\,I_{(F^2_{\C})^*}(x,x)=0\,\},\\
&\Hh_{half}:=\{x \in (\rad I_F)^*_{\C}\,|\,
	\langle a,x \rangle \neq 0,\,
	\langle b,x \rangle \neq 0,\,
	\mathrm{Im}\frac{\langle b,x \rangle}
	{\langle a,x \rangle}>0\,\}.	
\end{align}
These are complex manifolds and do not depend on the choice 
of an oriented basis $a,b$, but depend only on the orientation. 

We introduce group actions. 
Since the spaces $(F^2_{\C})^*$ and $(\rad I_F)^*_{\C}$ 
are $\C$-vector spaces, 
they have a natural $\C^*$-action. 
It induces $\C^*$-actions on $D^2$ and $\Hh_{half}$:
\begin{align}
\psi(\alpha):\ &D^2 \to D^2,\quad 
x \mapsto \alpha\cdot x \ (\forall \alpha \in \C^*),\\
\psi(\alpha):\ &\Hh_{half} \to \Hh_{half},\quad 
x \mapsto \alpha\cdot x \ (\forall \alpha \in \C^*).
\end{align}
The group $\widetilde{Aut}(R)$ acts on 
$(F^2_{\C})^*$ and $(\rad I_F)^*_{\C}$ by \EQUref{2.4action}. 
We could easily check that this action induces 
the actions 
\begin{align}
\varphi(g):\ &D^2 \to D^2\quad 
(\forall g \in \widetilde{Aut}(R)),\\
\varphi(g):\ &\Hh_{half} \to \Hh_{half}\quad 
(\forall g \in \widetilde{Aut}(R)).
\end{align}

We define $\C^*$-quotient spaces 
\begin{align}
&\PD:=D^2/\C^*,\\
&\PH:=\Hh_{half}/\C^*.
\end{align}
These spaces are complex manifolds and 
they have $\widetilde{Aut}(R)$-actions. 
The $\C^*$-bundles 
$p:D^2 \to \PD$, \hskip2mm$p:\Hh_{half} \to \PH$ 
are $\widetilde{Aut}(R)$-equivariant. 

We have a natural diagram:
\begin{equation}
\begin{CD}
D^2 @>>> \Hh_{half}\\
@V{p}VV @VV{p}V \\
\PD @>>> \PH
\end{CD}.\EQUlabel{2.4diagram}
\end{equation}
This diagram is a $\C^*$-equivariant morphism 
from $D^2 \to \PD$ to $\Hh_{half} \to \PH$. 
Thus this diagram is cartesian. 
Also this diagram is $\widetilde{Aut}(R)$-equivariant. 

\EQUnewpage

\subsection{A tensor on the domain $D^2$}
\subsubsection{Definition of a holomorphic metric on $D$}
First we define a holomorphic metric $I_D$ on $D$. 
Since $D$ is an open set of an affine space, 
there exists uniquely the holomorphic metric $I_{D}$ on $D$:
\begin{equation}
I_D \in \Gamma(D,\Omega^1_{D}\otimes_{\Oo_{D}}\Omega^1_{D})
\end{equation}
characterized by the property that on each $p \in D$, 
it gives 
\begin{equation}
I_{D,p}(\de,\de'):=
I_{(F^2_{\C})^*}(f_p(\de),f_p(\de'))
\quad
(\forall \de,\de' \in T_p D)
\end{equation}
where 
$
f_p:T_p D \stackrel{\sim}{\to}(F^2_{\C})^*
$ is a canonical isomorphism.

By this definition, we have 
\begin{align}
\varphi(g)^*I_D&=I_D\quad (\forall g \in \widetilde{Aut}(R)),
\EQUlabel{2.5action1}\\
\psi(\alpha)^*I_D&=\alpha^2I_D\quad (\forall \alpha \in \C^*).
\EQUlabel{2.5action2}
\end{align}
\EQUspace

\subsubsection{Definition of a tensor on $D^2$}
We define a tensor on $D^2$. 
We define 
$j_1^*I_D \in 
\Gamma(D^2,\Omega^1_{D^2}\otimes_{\Oo_{D^2}}\Omega^1_{D^2})$
as a pull-back of 
$I_D \in \Gamma(D,\Omega^1_D\otimes_{\Oo_D}\Omega^1_D)$
by the inclusion 
$j_1:D^2 \hookrightarrow D$. 
\begin{Prop}\EQUlabel{2.5prop}
There exists uniquely an element 
$I_{D^2}\in 
\Gamma(D^2,p^*\Omega^1_{\PD}\otimes_{\Oo_{D^2}}p^*\Omega^1_{\PD})$ 
such that its image by the morphism 
\begin{equation}
\Gamma(D^2,p^*\Omega^1_{\PD}\otimes_{\Oo_{D^2}}p^*\Omega^1_{\PD})
\to 
\Gamma(D^2,\Omega^1_{D^2}\otimes_{\Oo_{D^2}}\Omega^1_{D^2})
\end{equation}
is $j_1^*I_D$. 
The tensor $I_{D^2}$ gives a non-degenerate $\Oo_{D^2}$-symmetric 
bilinear form on $p^*\Theta_{\PD}$. 
\end{Prop}
\begin{proof}
For $x \in D^2$, the radical of $\C$-blinear form 
$(j_1^*I_D)_x:T_xD^2 \times T_xD^2 \to \C$ 
is just a kernel of 
$(p_*)_x:T_xD^2 \to T_{p(x)}\PD$. 
Thus $(j_1^*I_D)_x$ gives a non-degenerate $\C$-bilinear form 
on $T_xD^2/\mathrm{ker}(p_*)_x$. 
Thereby we have a result. 
\end{proof}
\EQUspace
By this proposition, we define its dual 
\begin{equation}
I^*_{D^2} \in 
\Gamma(D^2,p^*\Theta_{\PD}\otimes_{\Oo_{D^2}}p^*\Theta_{\PD}).
\end{equation}

By a uniqueness of $I_{D^2}$ and equations \EQUref{2.5action1}, 
\EQUref{2.5action2}, 
we have 
\begin{align}
\varphi(g)^*I^*_{D^2}&=I^*_{D^2}\quad (\forall g \in \widetilde{Aut}(R)),
\EQUlabel{2.5action3}\\
\psi(\alpha)^*I^*_{D^2}&=\alpha^{-2}I^*_{D^2}\quad (\forall \alpha \in \C^*).
\EQUlabel{2.5action4}
\end{align}

\EQUnewpage

\subsection{Definition of the analytic Weyl group quotient}
\subsubsection{Open subset of the domain}
For the open subsets of $D^2$ and $\PD$, 
we define $W^2$-quotient space and a tensor on it. 

We put 
\begin{align}
&\Dkai
	:=\{x \in D^2\,|\,
	\langle \alpha,x \rangle \neq 0\ 
	(\forall \alpha \in R)\,\},\\
&\Pp(\Dkai):=
	\Dkai/\C^*,
\end{align}
whose complement is a reflection 
hyperplanes. Then $p:\Dkai \to \Pp(\Dkai)$ is a 
$\C^*$-bundle and the diagram
\begin{equation}
\begin{CD}
\Dkai @>>> \Hh_{half} \\
@V{p}VV @VV{p}V \\
\Pp(\Dkai) @>>> \PH
\end{CD}
\end{equation}
is $\widetilde{Aut}(R)$-equivariant and cartesian. 
Since $j_2:\Dkai \to D^2$ is an open immersion, 
$j_2^*I^*_{D^2}$ is an element of 
$\Gamma(\Dkai,
p^*\Theta_{\Pp(\Dkai)}
\otimes_{\Oo_{\Dkai}}
p^*\Theta_{\Pp(\Dkai)})$. 
We put $I^*_{\Dkai}:=j_2^*I^*_{D^2}$. 
\begin{Prop}
The group $W^2$ acts on $\Dkai$ and $\Pp(\Dkai)$ 
properly discontinuous and fixed point free. 
\end{Prop}
\begin{proof}
Since $p:\Dkai \to \Pp(\Dkai)$ is $W^2$-equivariant, 
we should only prove this Proposition for 
$\Pp(\Dkai)$, which will be shown in Proposition \ref{5.3prop}. 
\end{proof}
\EQUspace
\subsubsection{Definition of an analytic Weyl group quotient}
Since $W^2$-action and $\C^*$-action are commutative, 
we have a $\C^*$-bundle 
$p:\Dkai/W^2 \to \Pp(\Dkai)/W^2$. 
In the following diagram:
\begin{equation}
\begin{CD}
\Dkai @>>{j_3}> \Dkai/W^2 @>>> \Hh_{half} \\
@V{p}VV @V{p}VV @VV{p}V \\
\Pp(\Dkai) @>>{j_3}> \Pp(\Dkai)/W^2 @>>> \PH
\end{CD},\EQUlabel{2.6diagram}
\end{equation}
all morphisms are $\widetilde{Aut}(R)$-equivariant 
and both squares are cartesian. 
\EQUspace

\subsubsection{A tensor on an analytic Weyl group quotient space}
We define a tensor on $\Dkai/W^2$. 
Since $I^*_{\Dkai}$ is $W^2$-invariant by \EQUref{2.5action3} and 
$j_3:\Dkai \to \Dkai/W^2$ is \'etale, 
there exists uniquely a tensor $I^*_{\Dkai/W^2} \in 
\Gamma(\Dkai/W^2,
p^*\Theta_{\Pp(\Dkai)/W^2}
\otimes_{\Oo_{\Dkai/W^2}}
p^*\Theta_{\Pp(\Dkai)/W^2})$ 
satisfying 
$j_3^*I^*_{\Dkai/W^2}=I^*_{\Dkai}$. 

By a uniqueness of $I^*_{\Dkai/W^2}$ and equations \EQUref{2.5action3}, 
\EQUref{2.5action4}, 
we have 
\begin{align}
\varphi(g)^*I^*_{\Dkai/W^2}&=I^*_{\Dkai/W^2}
\quad (\forall g \in \widetilde{Aut}(R)),\EQUlabel{2.6action1}\\
\psi(\alpha)^*I^*_{\Dkai/W^2}&=\alpha^{-2}I^*_{\Dkai/W^2}
\quad (\forall \alpha \in \C^*).\EQUlabel{2.6action2}
\end{align}

\EQUnewpage
\section{Algebraic Weyl group quotient spaces}

\subsection{Definition of Weyl group invariant rings}
\subsubsection{Definition of Weyl group invariants}
In this subsection, we define the Weyl group invariant ring for 
the oriented elliptic root system with  
$\Lambda_{\Z} \neq \emptyset$.

By the diagram \EQUref{2.4diagram}, we have the morphism:
\begin{align}
\pi:&\ D^2 \to \PH.
\end{align}

We define an $\Oo_{\PH}$-modules. For $k,m \in \Z$, we put 
\begin{align}
\Ss^W_{D^2,k,m}(U):=
	\{
	f \in \pi_*\Oo_{D^2}(U)|
	&\varphi(g)^* f=f\,(\forall g \in W^2),\,\\
	&\varphi(\rho(t))^*f=e^{-2\pi imt}f\,(\forall t \in \C),\nonumber\\
	&\psi(\alpha)^* f=\alpha^{-k} f\,(\forall \alpha \in \C^*)\,\},\nonumber
\end{align}
for an open set $U \subset \PH$. Here $\rho:\C \to K_{\C}$ 
is defined in Proposition \ref{2.3prop}. 

\subsubsection{Definition of Weyl group invariant rings}
We define the $\Oo_{\PH}$-graded algebras by 
\begin{equation}
\Ss^W_{D^2,*,*}:=\bigoplus_{k,m \in \Z}\Ss^W_{D^2,k,m},
\ \,
\Ss^W_{D^2,*,0}:=\bigoplus_{k \in \Z}\Ss^W_{D^2,k,0},
\ \,
\Ss^W_{D^2,0,*}:=\bigoplus_{m \in \Z}\Ss^W_{D^2,0,m}.
\end{equation}
\EQUspace

\subsubsection{Automorphism group action on the Weyl group invariant rings}
\EQUlabel{3.1subsection}
In this subsection, we define an $\widetilde{Aut}(R)$-action
and a $\C^*$-action on the ringed spaces 
$(\PH,\Ss^W_{D^2,*,*})$ etc. 

Since the group $\widetilde{Aut}(R)$ acts on $D^2$, 
we have 
\begin{equation}
\Oo_{D^2} \to \varphi(g)_*\Oo_{D^2},\quad
f \mapsto f \circ \varphi(g)
\end{equation}
for a local section $f$ of $\Oo_{D^2}$. 
Taking a direct image by $\pi$, we have 
\begin{equation}
\pi_*\Oo_{D^2} \to \pi_*\varphi(g)_*\Oo_{D^2}
\simeq \varphi(g)_*\pi_*\Oo_{D^2}.\EQUlabel{3.1morphism}
\end{equation}
A pair $\varphi(g):\PH \to \PH$ and the morphism 
\EQUref{3.1morphism} 
gives a morphism of a ringed space:
\begin{equation}
\varphi(g):(\PH,\pi_*\Oo_{D^2}) \to (\PH,\pi_*\Oo_{D^2})\ 
(g \in \widetilde{Aut}(R)).
\end{equation}
These give an $\widetilde{Aut}(R)$-action of a ringed space 
$(\PH,\pi_*\Oo_{D^2})$. 
Since we could restrict the morphism \EQUref{3.1morphism} 
to the subsheaves 
$\Ss^W_{D^2,*,*}$, $\Ss^W_{D^2,*,0}$, $\Ss^W_{D^2,0,*}$, 
and $\Ss^W_{D^2,0,0}$, 
we obtain $\widetilde{Aut}(R)$-actions on the ringed spaces 
\begin{align}
&(\PH,\Ss^W_{D^2,*,*}),\ (\PH,\Ss^W_{D^2,*,0}), \\
&(\PH,\Ss^W_{D^2,0,*}), \ (\PH,\Ss^W_{D^2,0,0})\nonumber
\end{align}
and we denote the action also by $\varphi$. 

In a same manner we could also define $\C^*$-actions 
on these ringed spaces and we denote the action also by $\psi$. 

By a parallel construction, we could also define 
$\widetilde{Aut}(R)$-actions and $\C^*$-actions 
on the ringed spaces 
\begin{align}
&(\PH,\Oo_{\PH}),\quad (\PH,\pi_*\Oo_{\Hh_{half}}), \\
&(\PH,\pi_*\Oo_{\Pp(\Dkai)/W^2}) ,\quad (\PH,\pi_*\Oo_{\Dkai/W^2}). 
\nonumber
\end{align}

\EQUnewpage
\subsection{Structure of Weyl group invariant rings}
In this subsection, we study the structure of the invariant ring. 
\subsubsection{A Chevalley's type theorem}
\begin{Thm}\EQUlabel{3.2theorem}
{\rm(
\cite{Chevalley3}, 
\cite{Chevalley4}, 
\cite{Chevalley6},
\cite{Chevalley2}, 
\cite{Chevalley1}, 
\cite{Chevalley5}
)}

The $\Oo_{\PH}$-graded algebra $\Ss^W_{D^2,0,*}$ is an 
$\Oo_{\PH}$-free algebra, i.e. 
\begin{equation}
\Ss^W_{D^2,0,*}=\Oo_{\PH}[s^1,\cdots,s^{n-1}]
\end{equation}
for $s^j \in \Gamma(\PH,\Ss^W_{D^2,0,*})$ with 
$c^1 \geq c^2 \geq \cdots \geq c^{n-1}>0$ 
and $n:=l+2$. 
We remark that $j$ of $s^j$ is a suffix.
\end{Thm}
As we see in \S \ref{5.6subsection}, we reduce its proof 
to the works  
\cite{Chevalley3}, 
\cite{Chevalley4}, 
\cite{Chevalley6},
\cite{Chevalley2}, 
\cite{Chevalley1}, 
\cite{Chevalley5}. 
\EQUspace

\begin{Cor}\EQUlabel{3.2cor}
\begin{enumerate}
\item We have a natural inclusion morphism 
$(\rad I_F)_{\C} \to \Gamma(\PH,\Ss^W_{D^2,-1,0})$. 
\item Take $a \in \rad I_F \setminus\{0\}$. 
By 1, we regard $a \in \Gamma(\PH,\Ss^W_{D^2,-1,0})$. 
We also have $a^{-1} \in \Gamma(\PH,\Ss^W_{D^2,1,0})$. 
Then we have isomoprhisms:
\begin{align}
\Ss^W_{D^2,*,*}&\simeq \C[a,a^{-1}]\otimes_{\C}\Ss^W_{D^2,0,*},\\
\Ss^W_{D^2,*,0}&\simeq \C[a,a^{-1}]\otimes_{\C}\Ss^W_{D^2,0,0}.
\end{align}
\item We have a canonical isomoprhism:
\begin{equation}
\Ss^W_{D^2,*,*}\simeq 
\Ss^W_{D^2,*,0}
\otimes_{\Ss^W_{D^2,0,0}}
\Ss^W_{D^2,0,*}.
\end{equation}
\item The $\Oo_{\PH}$-algebras 
$\Ss^W_{D^2,*,*}$, 
$\Ss^W_{D^2,*,0}$, 
$\Ss^W_{D^2,0,*}$ 
are of finite presentation over $\Oo_{\PH}$. 
\end{enumerate}
\end{Cor}
\begin{proof}
For 1, we first see that an element of 
$(\rad I_F)_{\C}$ is regarded as a holomorphic function 
on $\Hh_{half}$. 
By the morphism $D^2 \to \Hh_{half}$, we have 
$(\rad I_F)_{\C} \to \Gamma(\PH,\pi_*\Oo_{D^2})$. 
We could easily check that it gives 1. 
Since $a \in \Gamma(\PH,\Ss^W_{D^2,-1,0})$ is a nowhere vanishing holomorphic 
function, we have $a^{-1} \in \Gamma(\PH,\Ss^W_{D^2,-1,0})$. 
By $a^k \Ss^W_{D^2,k,m} \simeq \Ss^W_{D^2,0,m}$, we have 2. 
For 3, we show that a natural morphism 
$
\Ss^W_{D^2,*,0}
\otimes_{\Ss^W_{D^2,0,0}}
\Ss^W_{D^2,0,*}
\to 
\Ss^W_{D^2,*,*}
$ 
is an isomorphism. 
By 2, we have 
$
\Ss^W_{D^2,*,0}
\otimes_{\Ss^W_{D^2,0,0}}
\Ss^W_{D^2,0,*}
\simeq 
(\C[a,a^{-1}]\otimes_{\C}\Ss^W_{D^2,0,0})
\otimes_{\Ss^W_{D^2,0,0}}
\Ss^W_{D^2,0,*}
\simeq
\C[a,a^{-1}]\otimes_{\C}\Ss^W_{D^2,0,*}
\simeq
\Ss^W_{D^2,*,*}
$. 
The 4th assertion is a direct consequence of 
$\Ss^W_{D^2,0,0} \simeq \Oo_{\PH}$, 
2 and Theorem \ref{3.2theorem}. 
\end{proof}
\EQUnewpage
\subsection{Relation between Weyl group invariant rings}

We consider the spaces $D^2$, $\Hh_{half}$, 
$\Dkai$, $\Dkai/W^2$, ${\Pp(\Dkai)/W^2}$ 
as $\PH$-objects and we denote the structure morphisms by $\pi$.  
\begin{Prop}\EQUlabel{3.3prop}
\begin{enumerate}
\item By the morphisms $D^2 \to \Hh_{half} \to \PH$, 
we regard the $\Oo_{\PH}$-modules $\pi_*\Oo_{\Hh_{half}}$ and $\Oo_{\PH}$ 
as submodules of $\pi_*\Oo_{D^2}$. 
Then we have
	\begin{enumerate}
	\item  $\Ss^W_{D^2,0,0} = \Oo_{\PH}$,
	\item  $\Ss^W_{D^2,*,0} \subset \pi_*\Oo_{\Hh_{half}}$. 
	\end{enumerate}	
\item{}	By the morphisms $\Dkai \to \Dkai/W^2 \to \Pp(\Dkai)/W^2$, 
we regard the $\Oo_{\PH}$-modules 
$\pi_*\Oo_{\Dkai/W^2}$ and $\pi_*\Oo_{\Pp(\Dkai)/W^2}$ 
as submodules of $\pi_*\Oo_{\Dkai}$. 
We also regard $\Ss^W_{D^2,0,*}$ and $\Ss^W_{D^2,*,*}$ 
as  a submodules of $\pi_*\Oo_{\Dkai}$. 
Then we have
	\begin{enumerate}
	\item $\Ss^W_{D^2,0,*} \subset \pi_*\Oo_{\Pp(\Dkai)/W^2}$,
	\item $\Ss^W_{D^2,*,*} \subset \pi_*\Oo_{\Dkai/W^2}$.	
	\end{enumerate}
\item By the above inclusion morphisms, we have morphism of ringed spaces:
\begin{align}
&(\PH,\Oo_{\PH}) = (\PH,\Ss^W_{D^2,0,0}),\\
&(\PH,\pi_*\Oo_{\Hh_{half}}) \to (\PH,\Ss^W_{D^2,*,0}),\\
&(\PH,\pi_*\Oo_{\Pp(\Dkai)/W^2}) \to (\PH,\Ss^W_{D^2,0,*}),\\
&(\PH,\pi_*\Oo_{\Dkai/W^2}) \to (\PH,\Ss^W_{D^2,*,*}).
\end{align}
These morphisms are all $\widetilde{Aut}(R)$-equivariant 
and $\C^*$-equivariant. 	
\end{enumerate}
\end{Prop}
\begin{proof}
In Theorem \ref{3.2theorem}, $c^i>0$ for $1 \leq i \leq n-1$. 
Thus we have 1(a). 
By Corollary \ref{3.2cor}(2) and 
$a,a^{-1} \in \Gamma(\PH,\pi_*\Oo_{\Hh_{half}})$, 
we have 1(b). 
For 2, we could check the conditions directly. 
For 3, it is O.K. because our definitions of group actions 
are induced from $D^2$ and $\Dkai$. 
\end{proof}

\EQUnewpage

\subsection{Definition of Specan}
We remind the notion of $\Specan$ and see some properties. 

Let $\Cc$ be a category of $\Oo_{\PH}$-algebra 
which is of finite presentation, 
For an object $\Aa$ of the category $\Cc$, 
the analytic space $\Specan\,\Aa$ could be defined by \cite{Houzel}. 
It is characterized by a natural isomorphism:
\begin{equation*}
\mathrm{Hom}_{(An)/\PH}(X,\Specan\,\Aa)\simeq 
\mathrm{Hom}_{\Oo_{X}}(f^*\Aa,\Oo_X)
\end{equation*}
for an object $f:X \to \PH$ of the category $(An)/\PH$, 
which is a category of an analytic space with the structure 
morphism to $\PH$.

Since there exists 
a canonical isomorphism: 
$
\mathrm{Hom}_{\Oo_{X}}(f^*\Aa,\Oo_X)
\simeq 
\mathrm{Hom}_{\Oo_{\PH}}(\Aa,f_*\Oo_X)
$, 
we have
\begin{equation}
\mathrm{Hom}_{(An)/\PH}(X,\Specan\,\Aa)
\simeq 
\mathrm{Hom}_{\Oo_{\PH}}(\Aa,f_*\Oo_X).
\end{equation}

We denote the structure morphism 
$\Specan\,\Aa \to \PH$ by $\pi$. 

\EQUspace

\EQUspace

\EQUspace

\EQUnewpage
\subsection{Algebraic Weyl group quotient space}
\subsubsection{Definition of Algebraic Weyl group quotient space}
We define the Algebraic Weyl group quotient spaces by 
\begin{align}
D^2//W^2&:=\Specan\,\Ss^W_{D^2,*,*},\\
\PD//W^2&:=\Specan\,\Ss^W_{D^2,0,*}.
\end{align}
The groups $\widetilde{Aut}(R)$ and $\C^*$ 
act on the ringed spaces $(\PH,\Ss^W_{D^2,*,*})$ 
and $(\PH,\Ss^W_{D^2,0,*})$ as we see in \S \ref{3.1subsection}. 
Then a natural morphism $p:D^2//W^2 \to \PD//W^2$ has a 
structure of $\C^*$-bundle. 
Also the group 
$\widetilde{Aut}(R)$ acts as 
\begin{equation}
\begin{CD}
D^2//W^2 @>{\varphi(g)}>> D^2//W^2 \\
@V{p}VV @V{p}VV \\
\PD//W^2 @>{\varphi(g)}>> \PD//W^2 \\
@VVV @VVV \\
\PH @>{\varphi(g)}>> \PH
\end{CD},
\end{equation}
for $g \in \widetilde{Aut}(R)$. 
We remark that $\varphi(g):D^2//W^2 \to D^2//W^2$ 
and $\varphi(g):\PD//W^2 \to \PD//W^2$ are not 
$\PH$-morphisms. 

\EQUspace
\subsubsection{Relation to the analytic Weyl group quotient 
space}
\begin{Prop}\EQUlabel{3.5prop}
\begin{enumerate}
\item We have the canonical isomorphism
\begin{equation}
\Hh_{half} \simeq \Specan\,\Ss^W_{D^2,*,0}.
\end{equation}
\item We have a following natural diagram:
\begin{equation}
\begin{CD}
\Dkai/W^2
@>{j_4}>>
D^2//W^2
@>>>
\Hh_{half}\\
@V{p}VV @V{p}VV @V{p}VV \\
\Pp(\Dkai)/W^2
@>{j_4}>>
\PD//W^2
@>>>
\PH
\end{CD}.\EQUlabel{3.5diagram}
\end{equation}
All morphism are $\widetilde{Aut}(R)$-equivariant and 
both squares are cartesian.
\item The composite morphisms of row arrows coincide with 
morphisms in diagram \EQUref{2.6diagram}. 
\item In the diagram, the morphisms 
$j_4:\, \Dkai//W^2 \to D^2//W^2$ and $j_4:\, {\Pp(\Dkai)/W^2} \to \PD//W^2$ 
are open immersions and their images are open dense. 
\item An element $(\rad I_F)_{\C}$ gives a $\C^*$-equivariant 
holomorphic function on $\Hh_{half}$. 
Thus it also gives a $\C^*$-equivariant 
holomorphic function on $D^2//W^2$.
\end{enumerate}
\end{Prop}
\begin{proof}
We show 1. By Proposition \ref{3.3prop}(1)(b), 
we have $\Ss^W_{D^2,*,0} \subset \pi_*\Oo_{\Hh_{half}}$. 
Then we have a morphism $\Hh_{half} \to \Specan\,\Ss^W_{D^2,*,0}$. 
Since it is a $\C^*$-equivariant morphism of the $\C^*$-bundles 
$\Hh_{half} \to \PH$ and $\Specan\,\Ss^W_{D^2,*,0} \to \PH$, 
we see that it is an isomorphism. 
We show 2. By 
$
\Ss^W_{D^2,*,*}\simeq 
\Ss^W_{D^2,*,0}
\otimes_{\Ss^W_{D^2,0,0}}
\Ss^W_{D^2,0,*}
$
in Corollary \ref{3.2cor}(2) and 1, 
we see that the right hand side of the diagram is cartesian. 
We show 3. Since the morphism $X \to \Hh_{half}$ is determined 
by the induced mapping $(\rad I_F)_{\C} \to \Gamma(X,\Oo_X)$, 
we could check that the morphisms 
$\Dkai/W^2 \to \Hh_{half}$ in this diagram and \EQUref{2.6diagram} 
are same. 
Since the below row arrows of \EQUref{3.5diagram} 
are a $\C^*$-quotient of the above row arrows of \EQUref{3.5diagram}, 
we have a result.  
We show 4. 
The morphisms
$j_4:\,\Dkai//W^2 \to D^2//W^2$ and $j_4:\,{\Pp(\Dkai)/W^2} \to \PD//W^2$ 
could be defined by $\Ss^W_{D^2,0,*} \subset \pi_*\Oo_{\Pp(\Dkai)/W^2}$
and $\Ss^W_{D^2,*,*} \subset \pi_*\Oo_{\Dkai/W^2}$ in 
Proposition \ref{3.3prop}(2)(b). 
By Proposition \ref{5.7prop}, $j_4:{\Pp(\Dkai)/W^2} \to \PD//W^2$ 
is an open immersion and its image is open dense. 
Thus $j_4:\,\Dkai//W^2 \to D^2//W^2$
is also an open immersion and its image is open dense. 
The 5th assertion is a direct consequence of 2. 
\end{proof}
\EQUnewpage

\subsection{A tensor on the algebraic quotient space}
We assert an existence of a tensor on an algebraic quotient space 
$D^2//W^2$. 
\begin{Prop}\EQUlabel{3.6prop}
There exists uniquely an element 
\begin{equation}
I^*_{D^2//W^2}\in 
\Gamma(D^2//W^2,p^*\Theta_{\PD//W^2}
\otimes_{\Oo_{D^2//W^2}}p^*\Theta_{\PD//W^2})
\end{equation}
such that 
\begin{equation}
j_4^*I^*_{D^2//W^2}=I^*_{\Dkai/W^2}.
\end{equation}
\end{Prop}
We give a proof in \S \ref{5.8subsection}. 

By a uniqueness of $I^*_{D^2//W^2}$ and equations 
\EQUref{2.6action1}, \EQUref{2.6action2}, 
we have 
\begin{align}
\varphi(g)^*I^*_{D^2//W^2}&=I^*_{D^2//W^2}
\quad (\forall g \in \widetilde{Aut}(R)),\EQUlabel{3.6action1}\\
\psi(\alpha)^*I^*_{D^2//W^2}&=\alpha^{-2}I^*_{D^2//W^2}
\quad (\forall \alpha \in \C^*).\EQUlabel{3.6action2}
\end{align}

\EQUnewpage
\section{Results}
\subsection{Frobenius}
\subsubsection{Definition of Frobenius manifold}
\indent

We remind the definition of Frobenius manifold.
\begin{Def}\rm{(\cite[p.146, Def. 9.1]{Hertling})}\EQUlabel{3.333}
A Frobenius manifold is a tuple $(M,\circ,e,E,J)$ 
where $M$ is a complex manifold of dimension $\geq 1$ 
with holomorphic metric $J$ and multiplication $\circ$ 
on the tangent bundle, $e$ is a global unit field 
and $E$ is another global vector field, subject to the following 
conditions:
\begin{enumerate}
\item the metric is invariant under the multiplication, i.e., 
$J(X\circ Y,Z)=J(X,Y\circ Z)$
for local sections $X,Y,Z \in \Theta_M$,
\item (potentiality) the $(3,1)$-tensor 
$\nabla \circ $ is symmetric
(here, $\nabla$ is the Levi-Civita connection of the metric), i.e., 
$
\nabla_X(Y \circ Z)-Y \circ \nabla_X(Z)
-\nabla_Y(X \circ Z)+X \circ \nabla_Y(Z)-[X,Y]\circ Z=0,
$ for local sections $X,Y,Z \in \Theta_M$,
\item the metric $J$ is flat,
\item $e$ is a unit field and it is flat, i.e. $\nabla e=0$, 
\item the Euler field $E$ satisfies 
$Lie_{E}(\circ )=1 \cdot \circ$ and 
$Lie_{E}(J)=D \cdot J$ for some $D \in \C$.
\end{enumerate}
\end{Def}

\subsubsection{Definition of Intersection form}
\begin{Def}\rm{(\cite[p.191]{Dubrovin})}\EQUlabel{3.555}
For a Frobenius manifold $(M,\circ,e,E,J)$, 
we define an intersection form 
$I^*:\Omega^1_M \times \Omega^1_M \to \Oo_M$ 
by 
\begin{equation}\EQUlabel{2.9.3}
I^*(\omega_1,\omega_2)=
	J(E,J^*(\omega_1)\circ J^*(\omega_2))
\end{equation}
where $J^*:\Omega^1_{M} \to \Theta_{M}$
is the isomorphism induced by $J$. 
\end{Def}

\EQUnewpage
\subsection{Conformal Frobenius structure and the intersection form}
\subsubsection{Conformal Frobenius structure}
We define a notion of conformal Frobenius manifold. 
\begin{Def}
Let $L \stackrel{p}{\to}M$ be a $\C^*$-bundle. 
Let $(\circ,e,E,J)$ be a tuple where 
	\begin{align}
	\circ \in \Gamma(L,\ &
	p^*\Omega^1_{M}\otimes_{\Oo_{L}}
	p^*\Omega^1_{M}\otimes_{\Oo_{L}}
	p^*\Theta_{M}),	\\
	e \in \Gamma(L,\ &
	p^*\Theta_{M}), 	\\
	E \in \Gamma(L,\ &
	p^*\Theta_{M}), \\	
	J \in \Gamma(L,\ &
	p^*\Omega^1_{M}\otimes_{\Oo_{L}}
	p^*\Omega^1_{M}). 		
	\end{align}
We call a tuple $(\circ,e,E,J)$ a conformal Frobenius structure 
of $p:L \to M$ if 
it satisfies the following conditions:
\begin{enumerate}
	\item For $\C^*$-action, we have
	\begin{equation}
	\psi(\alpha)^* \circ=\circ,\ 
	\psi(\alpha)^* e =e,\ 
	\psi(\alpha)^* E =E,\ 
	\psi(\alpha)^* J =\alpha^2 J
	\end{equation}
	where $\psi(\alpha):L \to L$ is a $\C^*$-action for 
	$\alpha \in \C^*$. 
	\item $J$ is non-degenerate, i.e. $J$ gives 
	the isomorphism 
	\begin{equation}
	J:p^*\Theta_{M} \stackrel{\sim}{\to} p^*\Omega^1_{M}.
	\EQUlabel{4.2iso}
	\end{equation}
	\item There exists an open covering 
	$M=\cup_{\lambda \in \Lambda}M_{\lambda}$ 
	and sections 
	$\iota_{\lambda}:M_{\lambda} \to p^{-1}(M_{\lambda})$ 
	such that a tuple 
	$(M_{\lambda},\iota_{\lambda}^*\circ,\iota_{\lambda}^*e,
	\iota_{\lambda}^*E,\iota_{\lambda}^*J)$ 
	is a Frobenius manifold for each $\lambda \in \Lambda$ 
	under the identifications 
	\begin{equation}
	\iota_{\lambda}^*p^*\Omega^1_{M}\simeq \Omega^1_{M_{\lambda}},\quad
	\iota_{\lambda}^*p^*\Theta_{M}\simeq \Theta_{M_{\lambda}}.
	\end{equation}
\end{enumerate}
\end{Def}
\EQUspace
\subsubsection{A good section}
We deine a notion of a good section. 
\begin{Def}
Let $p:L \to M$ be a $\C^*$-bundle and $(\circ,e,E,J)$ be 
a conformal Frobenius structure. 
Let $U \subset M$ be an open set and 
$\iota:U \to p^{-1}(U)$ be a section. 
Then we call $(U,\iota)$ a good section if 
$(U,\iota^*\circ,\iota^*e,\iota^*E,\iota^*J)$ is a Frobenius manifold.
\end{Def}
\EQUspace

\subsubsection{Intersection form}
\begin{Def}
	We define the intersection form $I^* \in 
	\Gamma(L,p^*\Theta_M \otimes_{\Oo_L}p^*\Theta_M)$ 
	of a conformal Frobenius structure $(\circ,e,E,J)$ by 
	\begin{equation}
	I^*(\omega_1,\omega_2)=J(E,J^*(\omega_1)\circ J^*(\omega_2)),
	\end{equation}
	where we denote $J^*$ as an inverse of $J$ defined in 
	\EQUref{4.2iso}. 
\end{Def}
We remark that for a good section $(U,\iota)$, 
$\iota^*I^*$ gives an intersection form 
for a Frobenius manifold $(U,\iota^*\circ,\iota^*e,\iota^*E,\iota^*J)$. 

%
%
%

\EQUnewpage
\subsection{Sufficient condition}\EQUlabel{4.3subsection}
We prepare a sufficient condition for the existence 
of a conformal Frobenius structure with a given intersection form. 

We first fix some notations. 
For a trivial $\C^*$-bundle $p:L \to M$, 
the following data are equivalent:
\begin{enumerate}
	\item $L \simeq M \times \C^*$. 
	\item $\iota:M \to L$ satisfying $p \circ \iota=id.$
	\item $\C^*$-equivariant morphism $f:L \to \C^*$. 
\end{enumerate}
These are related as $M \times \C^* \simeq L,\ 
(x,\alpha) \mapsto \alpha\cdot\iota(x)$ 
and 
$L \simeq M \times \C^*,\ 
y \mapsto (p(x),f(y))$. 
We denote by $\iota(f):M \to L$ a section which corresponds to 
a $\C^*$-equivariant morphism $f:L \to \C^*$. 

We denote 
\begin{align}
V_M&:=
\Gamma(M,(\Theta_M)^{\otimes l} \otimes_{\Oo_M}(\Omega^1_M)^{\otimes m}),\\
V_L&:=
\Gamma(L,(p^*\Theta_M)^{\otimes l}
\otimes_{\Oo_L}(p^*\Omega^1_M)^{\otimes m}), \\
V_L^k&:=\{\omega \in V_L
\,|\,\psi(\alpha)^*\omega=\alpha^k \omega\}.
\end{align}
We put 
\begin{equation}
\Psi:V_M \to V_L,\quad
\omega \mapsto p^*\omega
\end{equation}
be a natural morphism. 
The following lemma is easy. 
\begin{Lemma}\EQUlabel{4.3lemma}
Let $p:L \to M$ be a trivial $\C^*$-bundle. 
\begin{enumerate}
\item $\Psi$ induces an isomorphism $V_M \simeq V_L^0$. 
\item Let $\iota(f)$ be a section which corresponds to 
a $\C^*$-equivariant morphism $f:L \to \C^*$. 
Then for any $k \in \Z$, the following isomorphisms 
\begin{align}
&V_M\simeq V_L^k,\quad 
\omega \mapsto f^kp^*\omega,\\
&V_L^k\simeq V_M,\quad 
\omega \mapsto \iota(f)^*\omega
\end{align}
are inverse to each other. 
\item Let $g:L \to \C^*$ be a $\C^*$-equivariant morphism. 
Then for $\omega \in V_L^k$, we have
\begin{equation}
\iota(g)^*\omega=(f/g)^k \iota(f)^* \omega,
\end{equation}
where we regard $f/g$ as a holomorphic function on $M$. 
\end{enumerate}
\end{Lemma}
By a correspondence given in Lemma \ref{4.3lemma}, we have 
a following sufficient condition for a existence 
of a conformal Frobenius manifold with an intersection form. 
\begin{Prop}\EQUlabel{4.3prop}
Let $p:L \to M$ be a trivial $\C^*$-bundle. 
Let $\iota(f)$ be a section which corresponds to 
a $\C^*$-equivariant morphism $f:L \to \C^*$. 
Let $I^*$ be an element 
of $\Gamma(L,p^*\Theta_M\otimes_{\Oo_L}p^*\Theta_M)$ 
with $\psi^*(\alpha)I^*=\alpha^{-2}I^*$ for $\alpha \in \C^*$. 
We assume that $(M,\circ,e,E,J)$ is a Frobenius manifold 
with an intersection form $\iota^*I^*$. 
Then $(p^*\circ, p^*e,p^*E,f^2p^*J)$ 
gives a conformal Frobenius structure with an intersection form 
$I^*$. 
\end{Prop}
\EQUnewpage
\subsection{Existence of a conformal Frobenius structure}
We assert the existence of a conformal Frobenius structure 
of $p:D^2//W^2 \to \PD/W^2$ with suitable conditions. 
\EQUspace

We first define a Euler operator on $\PD//W^2$. 
For a local section $f$ of $\Oo_{\PD//W^2}$, we put 
\begin{equation}
E_{\PD//W^2}f:=\frac{-1}{c^1}\frac{1}{2\pi\sqrt{-1}}
\lim_{t \to 0}\frac{\varphi(\rho(t))^*f-f}{t}.
\end{equation}
Then $E_{\PD//W^2}$ gives a holomorphic vector field 
on $\PD$. 
\EQUspace
For the oriented elliptic root system 
with $\Lambda_{\Z} \neq \emptyset$, 
the condition $c^1 >c^2$ is called ``codimension 1" in \cite{extendedII}, 
where $c^1$ and $c^2$ are defined in Theorem \ref{3.2theorem}. 

\begin{Thm}\EQUlabel{4.4theorem}
If the oriented elliptic root system 
with $\Lambda_{\Z} \neq \emptyset$ 
satisfies the condition of codimension 1, 
then there exists a conformal Frobenius structure 
$(\circ,e,E,J)$ of 
$p:D^2//W^2 \to \PD//W^2$ 
satisfying the following properties:
\begin{enumerate}
\item Its intersection form coincides with $I^*_{D^2//W^2}$. 
\item Its Euler field coincides with $E_{\PD//W^2}$ 
under the identification of Lemma \ref{4.3lemma}(1). 
\end{enumerate}
\end{Thm}
For a proof, we give it in \S \ref{5.9subsection}. 

\EQUnewpage

\subsection{Uniqueness of a conformal Frobenius structure}
We assert that a conformal Frobenius structure 
of $p:D^2//W^2 \to \PD//W^2$ is unique 
up to $\C^*$ multiplication under some conditions.
\begin{Thm}\EQUlabel{4.5theorem}
	Let $(\circ,e,E,J)$, $(\circ',e',E,J')$ be 
	conformal Frobenius structures 
	of $p:D^2//W^2 \to \PD/W^2$ 
	which satisfies conditions of 1 and 2
	of Theorem \ref{4.4theorem}. 
	We also assume that they have global good sections. 
	Then there exists $c \in \C^*$ such that \\
\begin{equation}
	(\circ',e',E,J')
	=(c^{-1}\circ,ce,E,c^{-1}J).
\end{equation}
\end{Thm}
For a proof, we give it in \S \ref{5.11subsection}. 

\EQUnewpage
\subsection{Good sections}
We classify good sections of a conformal Frobenius structure 
of Theorem \ref{4.4theorem}. 
\begin{Prop}\EQUlabel{4.6prop}
We fix a conformal Frobenius structure of $p:D^2//W^2 \to \PD//W^2$. 
\begin{enumerate}
\item For any $f \in (\rad I_F)_{\C}\setminus \{0\}$, 
it defines $f:D^2//W^2 \setminus \{f=0\} \to \C^*$
and the corresponding section 
$\iota(f):\PD//W^2 \setminus p(\{f=0\}) \to D^2//W^2\setminus \{f=0\}$. 
Then $\iota(f)$ is a good section. 
\item Let $(U,\iota)$ be a good section of 
this conformal Frobenius structure. 
Then the exists uniquely an element $f \in (\rad I_F)_{\C}\setminus \{0\}$ 
such that $\iota=\iota(f)|_U$. 
\end{enumerate}
\end{Prop}
For a proof, we give it in \S \ref{5.12subsection}. 

\EQUnewpage
\subsection{The automorphism group action}

We discuss a transformation of a conformal Frobenius structure 
by the action of $\widetilde{Aut}(R)$ on $D^2//W^2$.
\begin{Prop}
We fix a conformal Frobenius structure $(\circ,e,E,J)$ 
of $p:D^2//W^2 \to \PD//W^2$. 
Then there exists the group homomorphism
\begin{equation}
\chi:\widetilde{Aut}(R)\to \C^*
\end{equation}
such that  
\begin{equation}
(\varphi(g)^*\circ,\varphi(g)^*e,\varphi(g)^*E,\varphi(g)^*J)
=
(\chi(g)^{-1}\circ,\chi(g)e,E,\chi(g)^{-1}J)
\end{equation}
for any $g \in \widetilde{Aut}(R)$. 
\end{Prop}
\begin{proof}
We fix $g \in \widetilde{Aut}(R)$. 
We see that 
$(\varphi(g)^*\circ,\varphi(g)^*e,\varphi(g)^*E,\varphi(g)^*J)$ 
is a conformal Frobenius structure of 
$p:D^2//W^2 \to \PD//W^2$ with intersection form 
$\varphi(g)^*I^*_{D^2//W^2}$. 
By \EQUref{3.6action1}, we have $\varphi(g)^*I^*_{D^2//W^2}=I^*_{D^2//W^2}$. 
Since $E=E_{\PD//W^2}$ and $E_{\PD//W^2}$ is a center of 
$\widetilde{Aut}(R)$, 
we have $\varphi(g)^*E=E$. 
Then by a uniqueness of Theorem \ref{4.5theorem}, 
we see that 
\begin{equation}
(\varphi(g)^*\circ,\varphi(g)^*e,\varphi(g)^*E,\varphi(g)^*J)
=
(c(g)^{-1}\circ,c(g)e,E,c(g)^{-1}J). 
\end{equation}
for some $c(g) \in \C^*$. 
Since $\varphi(g)^*e=c(g)e$ and $e \neq 0$, 
$c(g)$ is uniquely determined and $c(g)$ satisfies 
$c(gg')=c(g)c(g')$. 
We put $\chi(g):=c(g)$. Then we have a result. 
\end{proof}
\EQUspace

\EQUnewpage
\section{Proofs}
\subsection{Domains for a signed marking}
Hereafter we fix an element $a \in \Lambda_{\Z}$. 
We take $b \in \rad I_F$ such that $a,b$ is an oriented basis. 
We put 
\begin{align}
&D_a:=\{x \in D\,|\,\langle a,x\rangle=1\,\},\\
&D^2_a:=D_a \cap D^2,\\
&\Hh_a:=\{x \in \Hh_{half}\,|\,\langle a,x\rangle=1\,\},\\
&F^a_{\C}:=\{x \in F^2_{\C}\,|\,I_{F^2_{\C}}(a,x)=0\,\},\\
&(F^a_{\C})^*:=\mathrm{Hom}_{\C}(F^a_{\C},\C),\\
&D^1_a:=\{x \in (F^a_{\C})^*\,|\,\langle a,x\rangle=1,\ 
\mathrm{Im}\langle b,x\rangle>0\,
\}.
\end{align}
The spaces $D_a$, $D^2_a$ and $D^1_a$ are complex manifolds
with natural $W^2$-action and $K_{\C}$-action. 
The space $D^1_a$ does not depend on the choice of $b \in \rad I_F$. 
\begin{Rem}\EQUlabel{5.1remark}
We remark that $F^a_{\C}$ gives a complexification of 
a hyperbolic extension defined in \cite{extendedII}. 
Thus $D^1_a$ corresponds to a space $\Ee$ defined in \cite{extendedII}.
\end{Rem}

We have a natural diagram:
\begin{equation}
\begin{CD}
D^1_a @<{f_2}<< D_a @<{f_1}<< D^2_a @>>> \Hh_a\\
@. @V{i}VV @V{i}VV @V{i}VV \\
@. D @<{j_1}<< D^2 @>>> \Hh_{half}\\
@. @. @V{p}VV @VV{p}V \\
@. @.\PD @>>> \PH
\end{CD}.\EQUlabel{5.1diagram}
\end{equation}
All three squares are cartesian. 
This diagram is $W^2$-equivariant and $K_{\C}$-equivariant. 
We also remark that the morphisms
\begin{equation}
f_2 \circ f_1:D^2_a \to D^1_a,\quad
p \circ i:D^2_a \to \PD
\end{equation}
are isomorphisms. 
\EQUnewpage
\subsection{A tensor on the domain}
We define holomorphic metrics on $D^2_a$ and $D^1_a$. 
We put
\begin{align}
&I_{D^2_a}:=i^*j_1^*I_D \in 
\Gamma(D^2_a,\Omega^1_{D^2_a}\otimes_{\Oo_{D^2_a}}\Omega^1_{D^2_a}),\\
&I_{D^1_a}:=(f_2 \circ f_1)^* I_{D^2_a}\in 
\Gamma(D^1_a,\Omega^1_{D^1_a}\otimes_{\Oo_{D^1_a}}\Omega^1_{D^1_a}).
\end{align}
By a proof of Proposition \ref{2.5prop}, 
$I_{D^2_a}$ is non-degenerate. 
Since $f_2 \circ f_1$ is an isomorphism, $I_{D^1_a}$ is also 
non-degenerate. 
We define their duals:
\begin{align}
&I^*_{D^2_a} \in 
\Gamma(D^2_a,\Theta_{D^2_a}\otimes_{\Oo_{D^2_a}}\Theta_{D^2_a}),\\
&I^*_{D^1_a} \in 
\Gamma(D^1_a,\Theta_{D^1_a}\otimes_{\Oo_{D^1_a}}\Theta_{D^1_a}).
\end{align}
We remark that $I^*_{D^1_a}$ is $W^2$-invariant. 
The following Proposition gives a characterization 
of the tensor $I^*_{D^1_a}$.
\begin{Prop}
For any $x \in D^1_a$, we have a canonical isomorphism
$T^*_xD^1_a \simeq F^a_{\C}/\C a$ and we have a commutative 
diagram:
\begin{equation}
\begin{CD}
	(I^*_{D^1_a})_x:\  & 
	T^*_xD^1_a & \times & 
	T^*_xD^1_a & @>>> 
	\C\\
	& @AAA @AAA & @|\\
	I_{F^a_{\C}/\C a}:\  & 
	F^a_{\C}/\C a & \times & 
	F^a_{\C}/\C a&
	@>>> \C
\end{CD},
\end{equation}
where $I_{F^a_{\C}/\C a}$ is a $\C$-bilinear mapping 
induced from $I_{F^2_{\C}}$. 
\end{Prop}
\begin{proof}
Take $y \in D^2_a$ such that $f_2\circ f_1(y)=x$.
Then we have
\begin{equation}
\begin{CD}
T_{x}D^1_a @<{(f_2)_*}<< T_{f_1(y)}D_a @<{(f_1)_*}<< T_{y}D^2_a\\
@. @V{i_*}VV @. \\
@. T_{i\circ f_1(y)}D @.
\end{CD}
\end{equation}
and a canonical isomorphisms:
\begin{equation}
\begin{CD}
T_{x}D^1_a @<(f_2)_*<< T_{f_1(y)}D_a @>i_*>> T_{i\circ f_1(y)}D\\
@A{\wr}AA @A{\wr}AA @A{\wr}AA \\
(F^a_{\C}/\C a)^* @<<< (F^2_{\C}/\C a)^* @>>> (F^2_{\C})^*
\end{CD}.
\end{equation}
From this, we have 
\begin{equation}
\mathrm{ker}[(f_2)_{*}:T_{f_1(y)}D_a \to T_{x}D^1_a]
=\rad (i^*I_D)_{f_1(y)}.
\end{equation}
Then $(i^*I_D)_{f_1(y)}$ induces a $\C$-bilinear form 
on $T_x D^1_a$. 
We see easily that it coincides with 
$I_{D^1_a}$ because $I_{D^1_a}$ could be written 
as $(f_2 \circ f_1)^*f_1^*i^*I_D$. 
\end{proof}
\begin{Rem}\EQUlabel{5.2remark}
By this Proposition, we see that a tensor $I^*_{D^1_a}$ coincides 
with a tensor $I^*_{\Ee}$ defined in \cite{extendedII} 
under the identification 
$D^1_a\simeq \Ee$. 
\end{Rem}
\begin{Rem}\EQUlabel{5.2remark2}
We also comment of a relation to the work \cite{Satakeauto}. 
In the diagram \EQUref{5.1diagram}, the group 
$\widetilde{Aut}(R)$ acts on $p:D^2 \to \PD$. 
But by $g \in \widetilde{Aut}(R)$, 
$D^2_a$ goes to $D^2_{g(a)}$. 
So we do not have a natural $\widetilde{Aut}(R)$ action on 
$D^2_a$ nor $D^1_a$. 
By the identification of $D^2_a$ and $D^1_a$ with 
$\PD$, we have an $\widetilde{Aut}(R)$-action on 
$D^2_a$ and $D^1_a$. 
Then the tensors $I^*_{D^2_a}$ and $I^*_{D^1_a}$ 
admit a conformal transformation by this action. 
This explains the reason of the appearance of 
a conformal transformation 
in the work of \cite{Satakeauto}.
\end{Rem}

\EQUnewpage

\subsection{Open subset of the domain 
and its analytic Weyl group quotient space}
We put
\newcommand{\Dakai}{\stackrel{\circ}{D^2_a}}
\newcommand{\Dakaione}{\stackrel{\circ}{D^1_a}}
\begin{align}
\Dakai:=\{x \in D^2_a\,|\,
\langle \alpha,x\rangle \neq 0\ (\forall \alpha \in R)\,\},\\
\Dakaione:=\{x \in D^1_a\,|\,
\langle \alpha,x\rangle \neq 0\ (\forall \alpha \in R)\,\}.
\end{align}
The group $W^2$ acts on $\Dakai$, $\Dakaione$ and we have 
a diagram:
\begin{equation}
\begin{CD}
\Dakaione @<<< \Dakai @>>> \Pp(\Dkai)\\
@VVV @VVV @VVV \\
D^1_a @<{f_2 \circ f_1}<< D^2_a @>{p \circ i}>> \Pp(D^2)
\end{CD}
\end{equation}
which is $W^2$-equivariant. 
\begin{Prop}\EQUlabel{5.3prop}
The group $W^2$ acts on $\Dakaione$, $\Dakai$, $\Pp(\Dkai)$ 
properly discontinuous and fixed point free. 
\end{Prop}
\begin{proof}
We need only to show it for $\Dakaione$. 
By Remark \EQUref{5.1remark}, 
we should prove for the domain $\Ee$ 
in \cite{extendedII}. 
It is proved in \cite{extendedII}. 
\end{proof}
\EQUnewpage

\subsection{A tensor on the analytic Weyl group quotient}
Since $j_3:\Dakai \to \Dakai/W^2$ and 
$j_3:\Dakaione \to \Dakaione/W^2$ are \'etale, 
there exist uniquely 
\begin{align}
I^*_{\Dakai/W^2}\in \Gamma(\Dakai/W^2,\Theta_{\Dakai/W^2}
\otimes_{\Oo_{\Dakai/W^2}}\Theta_{\Dakai/W^2}),\\
I^*_{\Dakaione/W^2}\in \Gamma(\Dakaione/W^2,\Theta_{\Dakaione/W^2}
\otimes_{\Oo_{\Dakaione/W^2}}\Theta_{\Dakaione/W^2})
\end{align}
such that 
\begin{equation}
I^*_{\Dakai}=j_3^*I_{\Dakai/W^2},\quad
I^*_{\Dakaione}=j_3^*I_{\Dakaione/W^2}.
\end{equation}
\EQUnewpage

\subsection{Weyl group invariant ring}
By the diagram \EQUref{5.1diagram}, 
we regard $\Dakai$, $\Dakaione$ and $\Hh_a$ as 
$\PH$-objects and we denote the structure morphisms 
by $\pi$. 
We put
\begin{align}
\Ss^W_{D^2_a,m}(U):=\{f \in \pi_*\Oo_{D^2_a}(U)\,|
&\,\varphi(g)^*f=f\,(\forall g \in W^2), \\
&\varphi(\rho(t))^*f=e^{-2\pi imt}f\,(\forall t \in \C)\,\},
\nonumber\\
\Ss^W_{D^1_a,m}(U):=\{f \in \pi_*\Oo_{D^1_a}(U)\,|
&\,\varphi(g)^*f=f\,(\forall g \in W^2), \\
&\varphi(\rho(t))^*f=e^{-2\pi imt}f\,(\forall t \in \C)\,\},
\nonumber
\\
\Ss^W_{\Hh_a}:=\pi_*\Oo_{\Hh_a}&
\end{align}
for an open set $U \subset \PH$. 
We define the $\Oo_{\PH}$-graded algebras by 
\begin{equation}
\Ss^W_{D^2_a,*}:=\oplus_{m \in \Z}\Ss^W_{D^2_a,m},\quad
\Ss^W_{D^1_a,*}:=\oplus_{m \in \Z}\Ss^W_{D^1_a,m}.
\end{equation}
\EQUnewpage

\subsection{Reduction of a proof of Theorem \ref{3.2theorem}
to the previous works}
\EQUlabel{5.6subsection}
Since $\Oo_{\PH}$-graded algebras 
$\Ss^W_{D^2_a,*}$, $\Ss^W_{D^1_a,*}$ and $\Ss^W_{D^2,0,*}$ 
are isomorphic, 
we need only to show Theorem \ref{3.2theorem} 
for $\Ss^W_{D^1_a,*}$, 
which is already shown in 
\cite{Chevalley3}, 
\cite{Chevalley4}, 
\cite{Chevalley6},
\cite{Chevalley2}, 
\cite{Chevalley1}, 
\cite{Chevalley5}.
\EQUnewpage

\subsection{Algebraic Weyl group quotient spaces}
We define the Weyl group quotient spaces by 
\begin{align}
D^2_a//W^2&:=\Specan\,\Ss^W_{D^2_a,*},\\
D^1_a//W^2&:=\Specan\,\Ss^W_{D^1_a,*}.
\end{align}
These spaces are related to the spaces which we defined in the previous 
sections.
\begin{Prop}\EQUlabel{5.7prop}
\begin{enumerate}
\item We have a canonical isomorphism:
\begin{equation}
\Hh_a \simeq \Specan\,\Ss^W_{\Hh_a}.
\end{equation}
\item We have a following natural diagram:
\begin{equation}
\begin{CD}
\Dakaione/W^2 @>{j_4}>> D^1_a//W^2 @>>> \Hh_a\\
@A{f_2 \circ f_1}AA @A{f_2 \circ f_1}AA @|\\
\Dakai/W^2 @>{j_4}>> D^2_a//W^2 @>>> \Hh_a\\
@V{i}VV  @V{i}VV  @V{i}VV  \\
\Dkai/W^2 @>{j_4}>> D^2//W^2 @>>> \Hh_{half}\\
@V{p}VV  @V{p}VV  @V{p}VV  \\
\Pp(\Dkai)/W^2 @>{j_4}>> \PD//W^2 @>>> \PH
\end{CD}.\EQUlabel{5.7diagram}
\end{equation}
These squares are all cartesian.
\item $j_4$ are all open immersion and the image is open dense. 
\item By Proposition \ref{3.5prop}(5), 
$a \in \Lambda_{\Z}$ defines a $\C^*$-equivariant morphism 
$a:D^2//W^2 \to \C^*$. 
By a correspondence in \S \ref{4.3subsection} 
it defines a section $\iota(a):\PD//W^2 \to D^2//W^2$. 
We have $\iota(a)=i \circ (p \circ i)^{-1}$. 
\end{enumerate}
\end{Prop}
\begin{proof}
A proof of 1,2 are parallel to Proposition \ref{3.5prop}.
For a proof of 3, we should only prove it for 
$j_4: \Dakaione/W^2 \to D^1_a//W^2$. 
By the identification in Remark \ref{5.1remark}, 
we reduce it to the corresponding result, 
which is shown in \cite{extendedII}. 
See also (3.12) and below of \cite{SatakeFrobenius}. 
The 4th assertion is a direct consequence of 2. 
\end{proof}

\EQUnewpage

\subsection{Proof of Proposition \ref{3.6prop}}
\EQUlabel{5.8subsection}
\begin{Prop}
There exists uniquely 
\begin{align}
I^*_{D^2_a//W^2}&\in \Gamma(D^2_a//W^2,\Theta_{D^2_a//W^2}
\otimes_{\Oo_{D^2_a//W^2}}\Theta_{D^2_a//W^2}),\\
I^*_{D^1_a//W^2}&\in \Gamma(D^1_a//W^2,\Theta_{D^1_a//W^2}
\otimes_{\Oo_{D^1_a//W^2}}\Theta_{D^1_a//W^2})
\end{align}
such that 
\begin{equation}
j_4^*I^*_{D^2_a//W^2}=I^*_{\Dakai/W^2},\quad
j_4^*I^*_{D^1_a//W^2}=I^*_{\Dakaione/W^2}.
\end{equation}
\end{Prop}
\begin{proof}

For the existence of $I^*_{D^1_a//W^2}$, we proved in \cite{SatakeFrobenius} 
under the identification in Remark \ref{5.1remark}, \ref{5.2remark}. 
We put 
\begin{equation}
I^*_{D^2_a//W^2}:=(f_2 \circ f_1)^*I^*_{D^1_a//W^2}.\EQUlabel{5.8equation}
\end{equation}
Then we have the result. 
\end{proof}

We define a tensor $I^*_{D^2//W^2}$ by 
\begin{equation}
I^*_{D^2//W^2}:=
a^{-2}p^*((p\circ i)^{-1})^*I^*_{D^2_a//W^2}.
\end{equation}
Then we have
\begin{align*}
j_4^*I^*_{D^2//W^2}
&=
j_4^*a^{-2}p^*((p\circ i)^{-1})^*I^*_{D^2_a//W^2}\\
&=
a^{-2}p^*((p\circ i)^{-1})^*j_4^*I^*_{D^2_a//W^2}\\
&=
a^{-2}p^*((p\circ i)^{-1})^*I^*_{\Dkai_a/W^2}\\
&=
a^{-2}p^*((p\circ i)^{-1})^*i^*I^*_{\Dkai/W^2}\\
&=
a^{-2}p^*(i \circ (p\circ i)^{-1})^*I^*_{\Dkai/W^2}\\
&=
a^{-2}p^*\iota(a)^*I^*_{\Dkai/W^2}\\
&=
I^*_{\Dkai/W^2},
\end{align*}
where we use the correspondence in Lemma \ref{4.3lemma} 
for the last equality. 

Since $\stackrel{\circ}{D^2}//W^2$ is open dense in $D^2//W^2$, 
the uniqueness of $I^*_{D^2//W^2}$ is apparent. 
Thus we obtain a proof of the Proposition \ref{3.6prop}. 
From the definition of $I^*_{D^2//W^2}$, we have the following 
easy consequence.

\begin{Cor}\EQUlabel{5.8cor}
By the identification $i^*p^*\Theta_{\PD//W^2}\simeq \Theta_{D^2_a//W^2}$, 
we have
\begin{equation}
i^*I^*_{D^2//W^2}=I^*_{D^2_a//W^2}.
\end{equation}
\end{Cor}

\EQUnewpage

\subsection{Proof of Theorem \ref{4.4theorem}}
\EQUlabel{5.9subsection}.
By the identification in Remark \ref{5.1remark} and \ref{5.2remark}, 
we have a following result by \cite{SatakeFrobenius}. 
\begin{Thm}\EQUlabel{5.9theorem}
The space $D^1_a//W^2$ has a structure of Frobenius manifold 
$(D^1_a//W^2,\circ,e,E,J)$
whose intersection form is $I^*_{D^2_a//W^2}$ 
and $E=E_{D^1_a//W^2}$, 
where we define $E_{D^1_a//W^2}$ by 
\begin{equation}
E_{D^1_a//W^2}f:=\frac{-1}{c^1}\frac{1}{2\pi\sqrt{-1}}
\lim_{t \to 0}\frac{\varphi(\rho(t))^*f-f}{t}
\end{equation}
for a local section $f \in \Oo_{D^1_a//W^2}$. 
We also assert that 
\begin{enumerate}
\item $J(e) \in \C^*d(b/a)$, 
\item $d(b/a)$ is flat for $I^*_{D^1_a//W^2}$
\end{enumerate} 
where we regard $b/a$ as a holomorphic function on 
$D^1_a//W^2$. 
\end{Thm}
Here we identify the operator $E_{D^1_a//W^2}$ 
with the operator in \cite{SatakeFrobenius} by the Remark \ref{2.3remark}.

By the isomorphims 
\begin{equation}
D^1_a//W^2 \stackrel{f_2 \circ f_1}{\longleftarrow}
D^2_a//W^2 \stackrel{p \circ i}{\longrightarrow}
\PD//W^2
\end{equation}
in the diagram \EQUref{5.7diagram}, 
we introduce the structure of Frobenius manifold on $\PD//W^2$. 
Its Euler field is $E_{\PD//W^2}$ because the diagram 
\EQUref{5.7diagram} is $K_{\C}$-equivariant. 
Its intersection form is 
\begin{align*}
&
\{(p \circ i)^{-1}\}^*(f_2 \circ f_1)^*I^*_{D^1_a//W^2}\\
=&
\{(p \circ i)^{-1}\}^*I^*_{D^2_a//W^2}\qquad
(\because \mbox{\EQUref{5.8equation}})\\
=&
\{(p \circ i)^{-1}\}^*i^*I^*_{D^2//W^2}\qquad
(\because \mbox{Corollary \ref{5.8cor}})\\
=&
\iota(a)^*I^*_{D^2//W^2}\qquad
(\because \mbox{Proposition \ref{5.7prop}(4)}).
\end{align*}
Then by Proposition \ref{4.3prop}, we have a proof of 
Theorem \ref{4.4theorem}. 
\EQUnewpage

\subsection{Deformation of holomorphic flat metric}
We give two Propositions. 
First is well-known and could be proved directly. 
\begin{Prop}\EQUlabel{5.10prop1}
Let $M$ be a complex manifold with dimension $n \geq 3$. 
Let $g$, $g'$ be flat holomorphic metrics on the complex manifold $M$
satisfying 
\begin{equation}
g'=\sigma^2 g
\end{equation}
for a nowhere-vanishing holomorphic function $\sigma$. 
Then $\sigma$ must be of the form
\begin{enumerate}
\item $\sigma^{-1}$ is constant,
\item $\sigma^{-1}=\sum_{ij}g_{ij}(x^i-c^i)(x^j-c^j)$ for a flat coordinate 
$x^1,\cdots,x^n$ with respect to $g$ and $c^i$ $(1 \leq i \leq n)$
\item $\sigma^{-1}=c-\sum_{i}b_{i}x^{i}$ 
for $c,b_i \in \C$ $(1 \leq i \leq n)$ satisfying 
$g^*(d\sigma^{-1},d\sigma^{-1})=0$ for a dual metric of $g$. 
\end{enumerate}
\end{Prop}
The following proposition is a Frobenius manifold version of 
the above proposition. 
\begin{Prop}\EQUlabel{5.10prop2}
Let $(M,\circ,e,E,J)$ be a Frobenius manifold 
with $\dim M \geq 3$, $J(e,e)=0$. 
Let $\sigma:M \to \C^*$ be a holomorphic function. \\
(1) $(M,\circ,e,E,J_{\sigma}:=\sigma^2J)$ is a Frobenius manifold 
if and only if $\sigma$ is homogeneous with respect to $E$ and 
$d(\sigma^{-1})=cJ(e)$ for some $c \in \C$. \\
(2) If $Lie_E J=D\cdot J$, then $Lie_E (\sigma^2 J)=(2-D)\sigma^2J$. 
\end{Prop}

\begin{proof}
We give a proof ``if'' part of (1). 
We check the conditions of Frobenius manifold. 
For a property of multiplication invariance of metric : 
\begin{equation}
J_{\sigma}(X \circ Y,Z)=J_{\sigma}(X,Y \circ Z), 
\end{equation}
it is O.K. by $J(X \circ Y,Z)=J(X,Y \circ Z)$. 

For a property of potentiality, we should prove that 
\begin{equation}\EQUlabel{7.2.1}
\nabla^{\sigma}_X (Y \circ Z)
-\nabla^{\sigma}_Y (X \circ Z)
+X \circ \nabla^{\sigma}_Y (Z)
-Y \circ \nabla^{\sigma}_X (Z)
-[X,Y]\circ Z=0
\end{equation}
for the Levi-Civita connection $\nabla^{\sigma}$ 
for a metric $J_{\sigma}$.  

Let $\nabla$ be the the Levi-Civita connection for a metric $J$. 
Since the relation of $\nabla^{\sigma}$ and $\nabla$ are known as 
\begin{equation}
\nabla^{\sigma}_X (Y)-\nabla_X(Y)=\alpha(X)Y+(X\wedge U)(Y)
\end{equation} 
where 
\begin{equation}\EQUlabel{7.12.3}
\alpha:=\frac{d\sigma}{\sigma},\ 
U:=J(\alpha),\ 
(X \wedge Y)(Z):=J(Y,Z)X-J(X,Z)Y.
\end{equation}
The equation \EQUref{7.2.1} holds 
if we replace $\nabla^{\sigma}$ with $\nabla$. 
Then we have
\begin{align*}
&\mathrm{LHS\ of\ }(\ref{7.2.1})\\
=&
 \alpha(X)(Y \circ Z)+(X \wedge U)(Y \circ Z)
-\alpha(Y)(X \circ Z)-(Y \wedge U)(X \circ Z)\\
&+X \circ (\alpha(Y)Z)+X \circ [(Y \wedge U)(Z)]
-Y \circ (\alpha(X)Z)-Y \circ [(X \wedge U)(Z)]\\
=&
(X \wedge U)(Y \circ Z)-(Y \wedge U)(X \circ Z)
+X \circ [(Y \wedge U)(Z)]-Y \circ [(X \wedge U)(Z)]\\
=&
J(U,Y \circ Z)X-J(U,X \circ Z)Y
-J(Y,Z)X \circ U+J(X,Z)Y \circ U.
\end{align*}
Here we used the property that a product $\circ$ is 
function linear, commutative. We also used that 
a metric $J$ is multiplication invariant. 

Since $U=-c\,\sigma\,e$, 
we have $J(U,Y \circ Z)X=J(Y,Z)X \circ U$ and 
$J(U,X \circ Z)Y=J(X,Z)Y \circ U$. Thus we see that 
\EQUref{7.2.1} holds. 

For a property of flatness for $J_{\sigma}$, 
we remind the relation between the curvature tensor of 
$J$ and $J_{\sigma}$. Let $R$ (resp. $R^{\sigma}$) be a curvature tensor 
of $J$ (resp. $J_{\sigma}$). Then 
\begin{equation}\EQUlabel{curvature}
R^{\sigma}(X,Y)=R(X,Y)-(B(X) \wedge Y+X \wedge B(Y)), 
\end{equation}
where $B(X):=-\alpha(X)U+\nabla_X(U)+\frac{1}{2}\alpha(U)X$, 
with $\alpha, U$ are defined in \EQUref{7.12.3}.
We have $R(X,Y)=0$ since $J$ is flat. 

In order to show that $B(X)=0$, we show that 
$\alpha(U)=0$ and $-\alpha(X)U+\nabla_X(U)=0$. 
The former part is as follows. 
$\alpha(U)=J(U,U)=c^2 \sigma^2J(e,e)$. 
Since $J(e,e)=0$, we have $\alpha(U)=0$. 
The latter part is as follows. 
$
-\alpha(X)U+\nabla_X(U)
=\frac{X(\sigma^{-1})}{\sigma^{-1}}U+\nabla_X(U)
=\sigma\nabla_X(\sigma^{-1}U)
=-c\,\sigma\nabla_X (e)
=0
$. 
Thus we have $B(X)=0$. Thereby $R^{\sigma}(X,Y)=0$. 

For a property $\nabla^{\sigma} e=0$, 
we use the equation 
\begin{equation}\EQUlabel{nabla}
\nabla^{\sigma}_X(e)-\nabla_X (e)=
\alpha(X)e+(X \wedge U)(e)
=J(U,X)e+J(U,e)X-J(X,e)U.
\end{equation}
First we have $\nabla_X (e)=0$. 
Since $U=-c\,\sigma\,e$, we have 
$J(U,X)e-J(X,e)U=0$. 
Since $J(U,e)=-c\,\sigma(e,e)=0$, we have 
$\nabla^{\sigma} e=0$. 

For a property of homogeneity conditions 
$Lie_E(J_a)=J_a$, we obtain it because 
$Lie_E(J)=J$ and $Lie_E(\sigma)=0$. 

We give a proof of ``only if'' part of (1). 
We assume that $(M,\circ,e,E,J_{\sigma})$ is a Frobenius manifold. 

We first show that $U=fe$ for some $f \in \Gamma(M,\Oo_M)$.
By the conditions $\nabla_X e=\nabla^{\sigma}_Xe=0$ 
and \EQUref{nabla}, we have 
\begin{equation}\EQUlabel{nabla2}
J(U,X)e+J(U,e)X-J(X,e)U=0. 
\end{equation}
We substitute $X=e$, we have $2J(U,e)e=0$. 
Since $e \neq 0$, we have $J(U,e)=0$. 
By \EQUref{nabla2}, we have $J(U,X)e=J(X,e)U$. 
Since $J$ is non-degenerate and flat, 
we take $X$ such that $J(X,e)$ is non-zero constant. 
Then we have $U=fe$ for some $f \in \Gamma(M,\Oo_M)$. 

We show that $U=-c\,\sigma\,e$ for some $c \in \C$. 
By \EQUref{curvature}, we have 
\begin{equation}\EQUlabel{curvature2}
B(X)\wedge Y+X \wedge B(Y)=0. 
\end{equation}
We could easily check that if $\dim M \geq 3$ then 
\EQUref{curvature2} implies that $B(X)=0$. 
Since $U=fe$, $\alpha(U)=J(U,U)=f^2J(e,e)=0$. 
Thus we have $-\alpha(X)U+\nabla_X(U)=0$. 
Since 
$
-\alpha(X)U+\nabla_X(U)
=\sigma\nabla_X(\sigma^{-1}U)
=\sigma\nabla_X(\sigma^{-1}fe)
=\sigma X(\sigma^{-1}f)e
$, 
we have $\sigma^{-1}f=-c$ for some $c \in \C$. 
Thus we have a result. 

For a proof of (2), we see that $Lie_E(J(e))=(D-1)J(e)$. 
Then we have $Lie_E(\sigma)=(1-D)\sigma$ and 
$Lie_E(\sigma^2J)=(2-2D+D)(\sigma^2J)=(2-D)(\sigma^2J)$. 
\end{proof}
\begin{Cor}\EQUlabel{5.10cor}
For a conformal Frobenius manifold $(\circ,e,E,J)$ 
of a $\C^*$-bundle $p:L \to M$, 
we assume that $\iota(f)$ is a good section 
for a $\C^*$-equivariant morphism 
$f:L \to \C^*$. 
For $g:L \to \C^*$: a $\C^*$-equivariant morphism, 
$\iota(g)$ is a good section if and only if 
\begin{equation}
d((f/g)^{-1})=c (\iota(f)^*J)(\iota(f)^*e)
\end{equation}
for some $c \in \C$. 
\end{Cor}
\begin{proof}
Since $\iota(g)^*J=(f/g)^2\iota(f)^*J$ by Lemma \ref{4.3lemma}, 
we obtain the result. 
\end{proof}
\EQUnewpage

\subsection{Proof of Theorem \ref{4.5theorem}}
\EQUlabel{5.11subsection}
We may assume that $(\circ,e,E,J)$ be a conformal Frobenius 
structure constructed in \S \ref{5.9subsection}. 
Then it has a good section $\iota(a)$. 
Let $h:D^2//W^2 \to \C^*$ be a $\C^*$-equivariant mapping whose 
section $\iota(h)$ is a good section of $(\circ',e',E,J')$. 
\begin{Lemma}\EQUlabel{5.11lemma1}
By the embedding 
$(\rad I_F)_{\C} \hookrightarrow \Gamma(D^2//W^2,\Oo_{D^2//W^2})$ 
in Proposition \ref{3.5prop}(4), $h$ is the image of this embedding. 
\end{Lemma}
\begin{proof}
On $\Pp(\Dkai)//W^2$, 
$\iota(h)^*I^*_{D^2//W^2}$ and $\iota(a)^*I^*_{D^2//W^2}$ 
are non-degenerate dual metric. 
Also they are flat because the intersection form 
of the Frobenius manifold is flat (\cite[p.191]{Dubrovin}). 

Since 
\begin{equation}
\iota(h)^* I^*_{D^2//W^2}
=(a/h)^{- 2}
\iota(a)^*I^*_{D^2//W^2}
\end{equation}
by Lemma \ref{4.3lemma}, 
we apply Proposition \ref{5.10prop1} and we see $(h/a)$ must be constant 
or suitable quardic or a suitable affine function with respect to 
flat coordinates for $\iota(a)^*I^*_{D^2//W^2}$. 
Also $h/a$ must be an element of $\Ss^W_{D^2,0,*}$ 
because the tensors $\iota(h)^*I^*_{D^2//W^2}$ and 
$\iota(a)^*I^*_{D^2//W^2}$ are homogeneous. 

By the isomorphism $\PD//W^2 \simeq D^1_a//W^2$ of the diagram 
\EQUref{5.7diagram}, $\iota(a)^*I^*_{D^2//W^2}$ corresponds to 
$I^*_{D^1//W^2}$ on $D^1//W^2$. 
Since flat coordinates for $I^*_{D^1//W^2}$ is an affine 
coordinate of $D^1$, we see that 
$h/a$ must be of the form $x(b/a)a+y$ for 
$b \in (\rad I_F)_{\C}\setminus \C a$ and $x,y \in \C$. 
Thus we have a result.
\end{proof}

We remind the fact. 
For a Frobenius manifold $(M,\circ,e,E,J)$ with intersection 
form $I^*$, we put 
\
\begin{equation}
\Ff_M(v):=
\{\omega \in \Omega^1_{M}\,|\,
Lie_{v}\omega=0\,\}
\end{equation}
for $v \in \Gamma(M,\Theta_{M})$. 
In \cite[Proposition 5.1]{SatakeFrobenius}, we have shown the fact that 
\begin{align}
&e^2I^*(\omega_1,\omega_2)=0,\EQUlabel{5.11equation1}\\
&eI^*(\omega_1,\omega_2)=J^*(\omega_1,\omega_2)
\EQUlabel{5.11equation2}
\end{align}
for $\omega_1,\omega_2 \in \Ff_M(e)$. 
We have a following Lemma. 
\begin{Lemma}\EQUlabel{5.11lemma2}
For a conformal Frobenius structure 
$(\circ',e',E,J')$, we have
\begin{equation}
(\iota(a)^*e')^2
[\iota(a)^*I^*_{D^2//W^2}(\omega_1,\omega_2)]
=0
\end{equation}
for $\omega_1,\omega_2 \in \Ff_{\PD//W^2}(\iota(a)^*e')$. 
\end{Lemma}
\begin{proof}
Since 
$(\PD//W^2,\iota(h)^*\circ',\iota(h)^*e',\iota(h)^*E,\iota(h)^*J')$ 
is a Frobenius manifold with intersection form 
$\iota(h)^*I^*_{D^2//W^2}$, 
we have 
\begin{equation}
(\iota(h)^*e')^2
[\iota(h)^*I^*_{D^2//W^2}(\eta_1,\eta_2)]=0
\end{equation}
for $\eta_1,\eta_2 \in \Ff_{\PD//W^2}(\iota(h)^*e')$ 
by \EQUref{5.11equation1}. 

Since 
$\iota(a)^*e'=\iota(h)^*e'$ and 
$\iota(a)^*I^*_{D^2//W^2}=(h/a)^{-2}\iota(h)^*I^*_{D^2//W^2}$
by Lemma \ref{4.3lemma}, 
we have
\begin{align*}
&(\iota(a)^*e')^2
[\iota(a)^*I^*_{D^2//W^2}(\omega_1,\omega_2)]\\
=&
(\iota(h)^*e')^2
[(h/a)^{-2}\iota(h)^*I^*_{D^2//W^2}(\omega_1,\omega_2)]\\
=&
(h/a)^{-2}
(\iota(h)^*e')^2
[\iota(h)^*I^*_{D^2//W^2}(\omega_1,\omega_2)]\\
=&0
\end{align*}
for $\omega_1,\omega_2 \in 
\Ff_{\PD//W^2}(\iota(a)^*e')$, where 
we used the fact that $(\iota(h)^*e')(h/a)=0$ 
because its degree is negative with respect to the Euler field 
$\iota(h)^*E$ by Lemma \ref{5.11lemma1}. 
\end{proof}
\begin{Lemma}\EQUlabel{5.11lemma3}
There exists $c \in \C^*$ such that 
$e'=ce$, $J'=c^{-1}J$. 
\end{Lemma}
\begin{proof}
Since $(\PD//W^2,\iota(a)^*\circ,\iota(a)^*e,\iota(a)^*E,\iota(a)^*J)$ 
is a Frobenius manifold with intersection form 
$\iota(a)^*I^*_{D^2//W^2}$, 
we have 
\begin{equation}
(\iota(a)^*e)^2
[\iota(a)^*I^*_{D^2//W^2}(\omega_1,\omega_2)]
=0\EQUlabel{5.11equation3}
\end{equation}
for $\omega_1,\omega_2 \in \Ff_{\PD//W^2}(\iota(a)^*e)$ 
by \EQUref{5.11equation1}. 

As we have seen in \cite{extendedII}, 
a non-singular vector field $v \in D^1//W^2$ of degree $-1$ 
satisfying 
\begin{equation}
v^2I^*_{D^1_a//W^2}(\omega_1,\omega_2)=0
\end{equation}
for $\omega_1,\omega_2 \in \Ff_{D^1_a//W^2}(v)$ 
is unique up to $\C^*$-multiplication, see also 
\cite[(4.13)]{SatakeFrobenius}. 
By Lemma \ref{5.11lemma2} and \EQUref{5.11equation3}, 
we have $\iota(a)^* e'=c\iota(a)^* e$ for some $c \in \C^*$. 

By \EQUref{5.11equation2}, we have 
\begin{align*}
(\iota(h)^*e')(\iota(h)^*I^*_{D^2//W^2})(\omega_1,\omega_2)
&=(\iota(h)^*J')^*(\omega_1,\omega_2)
\mbox{ for }
\omega_1,\omega_2\in \Ff_{\PD//W^2}(\iota(h)^*e'),
\\
(\iota(a)^*e)(\iota(a)^*I^*_{D^2//W^2})(\omega_1,\omega_2)
&=(\iota(a)^*J)^*(\omega_1,\omega_2)
\mbox{ for }
\omega_1,\omega_2 \in \Ff_{\PD//W^2}(\iota(a)^*e).
\end{align*}
Since $\iota(h)^*e'=\iota(a)^*e'=c\iota(a)^*e$, 
we have $\Ff(\iota(h)^*e')=\Ff(\iota(a)^*e)$. 
Then we have 
\begin{align*}
&(\iota(a)^*J)^*(\omega_1,\omega_2)\\
=&(\iota(a)^*e)(\iota(a)^*I^*_{D^2//W^2})(\omega_1,\omega_2)\\
=&c^{-1}(\iota(h)^*e')[(h^2/a^2)\iota(h)^*I^*_{D^2//W^2}](\omega_1,\omega_2)\\
=&c^{-1}(h^2/a^2)(\iota(h)^*J')^*(\omega_1,\omega_2)\\
=&c^{-1}(\iota(a)^*J')^*(\omega_1,\omega_2)
\end{align*}
for $\omega_1,\omega_2 \in \Ff(\iota(a)^*e)$. 
Since $\Ff(\iota(a)^*e)$ is sufficient to determine 
the tensor (\cite{extendedII}), 
we have $(\iota(a)^*J)^*=c^{-1}(\iota(a)^*J')^*$. 
By the correspondence of Lemma \ref{4.3lemma}, 
we have $e'=ce$ and $J'=c^{-1}J$ . 
\end{proof}
\begin{Lemma}
$\iota(a)$ is a good section of $(\circ',e',E,J')$. 
\end{Lemma}
\begin{proof}
By Lemma \ref{5.11lemma1}, $h \in (\rad I_F)_{\C}$. 
By Theorem \ref{5.9theorem}, 
$d(b/a)=const. (\iota(a)^*J)(\iota(a)^*e)$. 
Thus we have $d(h/a)=c_1(\iota(a)^*J)(\iota(a)^*e)$ for some $c_1 \in \C$. 
Then we have
\begin{align*}
d[(h/a)^{-1}]
&=
-(h/a)^{-2}d(h/a)\\
&=
-(h/a)^{-2}[c_1 (\iota(a)^*J)(\iota(a)^*e)]
\\
&=
-c_1(\iota(h)^*J)(\iota(h)^*e)\\
&=
-c_1(\iota(h)^*J')(\iota(h)^*e')
\quad (\because 
\mbox{ Lemma \ref{5.11lemma3}}).
\end{align*}
Since $\iota(h)$ is a good section, 
we apply the Corollary \ref{5.10cor}, then 
we have a result. 
\end{proof}

We give a proof of Theorem \ref{4.5theorem}.
Then $(\PD//W^2,\iota(a)^*\circ,\iota(a)^*e,\iota(a)^*E,\iota(a)^*J)$ 
and  $(\PD//W^2,\iota(a)^*\circ',\iota(a)^*e',\iota(a)^*E,\iota(a)^*J')$ 
are Frobenius manifolds whose intersection forms are 
$\iota(a)^*I^*_{D^2//W^2}$. 
Since we proved in \cite{SatakeFrobenius} 
that Frobenius manifold structure 
whose intersection form is $I^*_{D^1_a//W^2}$ and Euler field 
$E_{D^1_a//W^2}$ is unique up to $\C^*$ on $D^1_a//W^2$, 
we also have  
\begin{equation}
(\iota(a)^*\circ,\iota(a)^*e,\iota(a)^*E,\iota(a)^*J)
=
(c^{-1}\iota(a)^*\circ',c\iota(a)^*e',\iota(a)^*E,c^{-1}\iota(a)^*J')
\end{equation}
for some $c \in \C^*$. 
By the correspondence in Lemma \ref{4.3lemma}, we obtain the result. 
\EQUnewpage

\subsection{Proof of Proposition \ref{4.6prop}}
\EQUlabel{5.12subsection}
By a uniqueness of a conformal Frobenius structure, 
we may assume that $(\circ,e,E,J)$ is a conformal 
Frobenius structure constructed in \S \ref{5.9subsection}. 
Then $(\iota(a)^*J)(\iota(a)^*e)$ equals to $d(b/a)$ 
up to a constant multiple by Theorem \ref{5.9theorem}. 
Since $\iota(a)$ is a good section, 
we apply Corollary \ref{5.10cor} 
and we have a result. 
\EQUnewpage

\end{document}